\input amstex
\documentstyle{amsppt}
%
\input pictex
\def\begingraph<#1,#2>[#3]{%
  \beginpicture
  \setcoordinatesystem units <#1, #2>
  \setlinear
  \linethickness=#3
  }
\def\endgraph{%
  \endpicture
  }
\def\beginsmallgraph{\begingraph<0.5cm,0.5cm>[0.35pt]}

\def\vertex at (#1, #2){\put{$\bullet$} at #1 #2}
\catcode`!=11 
\def\edge from (#1, #2) to (#3, #4){%
  \!xloc=\!M{#1}\!xunit  \!xxloc=\!M{#3}\!xunit%
  \!yloc=\!M{#2}\!yunit  \!yyloc=\!M{#4}\!yunit%
  \!dxpos=\!xxloc  \advance\!dxpos by -\!xloc
  \!dypos=\!yyloc  \advance\!dypos by -\!yloc
  \ifdim\!dypos=\!zpt
    \putrule from #1 #2 to #3 #4
  \else
    \ifdim\!dxpos=\!zpt
      \putrule from #1 #2 to #3 #4
    \else
      \plot #1 #2 #3 #4 /
    \fi
  \fi
  }
\catcode`!=12 
\def\label #1 at (#2, #3){\put {#1} at #2 #3}

\magnification=\magstep1
\hsize=6.5 true in
\vsize=8.6 true in

\def\A{{\Cal A}}

\def\CC{{\Cal C}}
\def\L{{\Cal L}}
\def\Y{{\Cal Y}}
\def\RR{{\bold R}}
\def\S{{\Cal S}}
\def\W{{\Cal W}}
\def\P{{\Pi}}
\def\Z{{\Bbb Z}}
\def\R{{\Bbb R}}
\def\C{{\Bbb C}}
\def\CP{{\Bbb{CP}}}

\def\a{{\alpha}}
\def\b{{\beta}}

\def\l{{\lambda}}
\def\s{{\sigma}}

\def\im{\operatorname{im}}

\def\id{\operatorname{id}}
\def\pr{\operatorname{pr}}
\def\Aut{\operatorname{Aut}}
\def\Out{\operatorname{Out}}

\def\GL{\operatorname{GL}}

\def\Real{\operatorname{Re}}
\def\Imag{\operatorname{Im}}
\def\Hom{\operatorname{Hom}}
\def\conj{\operatorname{conj}}

\topmatter

\title{The Braid Monodromy of Plane Algebraic Curves
and Hyperplane Arrangements}
\endtitle

\rightheadtext{Braid Monodromy of Curves and Arrangements}
\leftheadtext{Daniel C.~Cohen and Alexander I.~Suciu}

\author Daniel C.~Cohen and Alexander I.~Suciu
\endauthor

\address{Daniel C. Cohen, Department of Mathematics,
Louisiana State Univeristy, Baton Rouge, LA  70808  USA}
\endaddress
\email{cohen\@math.lsu.edu}
\endemail
\address{Alexander I. Suciu, Department of Mathematics,
Northeastern Univeristy, Boston, MA  02115  USA}
\endaddress
\email{alexsuciu\@neu.edu}
\endemail

\thanks{The second author was partially supported by N.S.F.~grant
DMS--9504833, and an RSDF~grant from Northeastern University.}
\endthanks

\endsubjclass

\endkeywords

\abstract
To a plane algebraic curve of degree $n$, Moishezon associated
a braid monodromy homomorphism from a finitely generated free group
to Artin's braid group $B_n$.  Using Hansen's polynomial covering
space theory, we give a new interpretation of this construction.
Next, we provide an explicit description of the braid monodromy of an 
arrangement of complex affine hyperplanes, 
by means of an associated ``braided wiring diagram.''  The ensuing
presentation of the fundamental group of the complement is 
shown to be Tietze-I
equivalent to the Randell-Arvola presentation.   Work of Libgober 
then implies that the complement of a line arrangement is homotopy 
equivalent to the 2-complex modeled on either of these presentations.  
Finally, we prove that the braid monodromy of a line arrangement 
determines the intersection lattice.  Examples of Falk then show 
that the braid monodromy carries more information than the 
group of the complement, thereby answering a question of Libgober.
\endabstract

\endtopmatter

\document

\head 1. Introduction
\endhead

\subhead 1.1
\endsubhead
Let $\CC$ be an algebraic curve in $\C^2$.  In the 1930's,
Zariski commissioned van~Kampen to compute the fundamental group
of the complement, $\pi_1(\C^2 \setminus \CC)$.  The algorithm
for doing this was developed in \cite{vK}.  Refinements of
van~Kampen's algorithm were given by Chisini in the 50's,
and Ch\'eniot, Abelson, and Chang in the 70's.  In the early 80's,
Moishezon \cite{Mo} introduced the notion of braid monodromy,
which he used to recover van~Kampen's presentation.  Finally,
Libgober \cite{L1} showed that the $2$-complex associated to
the braid monodromy presentation is homotopy equivalent to 
$\C^2 \setminus \CC$.

Let $\A$ be an arrangement of hyperplanes in $\C^\ell$.  In the
early 80's, Randell \cite{R1} found an algorithm for computing the
fundamental group of the complement, $\pi_1(\C^\ell \setminus \A)$,
when $\A$ is the complexification of a real arrangement.
Salvetti~\cite{S1} subsequently found a regular cell complex that is
a deformation retract of the complement of such an arrangement.
When $\ell=2$, Falk \cite{Fa} proved that the $2$-complex associated
to the Randell presentation is homotopy equivalent to $\C^2 \setminus
\A$ by showing that it is homotopy equivalent to Salvetti's complex.
The braid monodromy of a complexified real arrangement was determined
by Salvetti~\cite{S2}, Hironaka~\cite{Hir}, and Cordovil and
Fachada~\cite{CF}, \cite{Cor}. An algorithm for computing the
fundamental group of an arbitrary complex arrangement was found by
Arvola \cite{Ar} (see also Orlik and Terao \cite{OT}, and see
Dung and Ha \cite{DH} for another method).

In this paper, we present a unified view of these two subjects,
extending several of the aforementioned results.  In particular,
we give in~5.3 an algorithm for finding the (pure) braid monodromy
of an arbitrary arrangement $\A$ of complex lines in $\C^2$.
Furthermore, we show in Theorem~6.4 that the corresponding
presentation of $\pi_1(\C^2 \setminus \A)$ is equivalent to the
Randell-Arvola presentation.  We also strengthen
Falk's result, by showing that the $2$-complex modeled on the
Arvola presentation is homotopy equivalent to $\C^2 \setminus \A$.

The determination of the braid monodromy of an arrangement $\A$
is facilitated by use of a {\it braided wiring diagram}
associated to $\A$, a natural generalization of
a combinatorial notion of Goodman \cite{Go}.  For a real arrangement,
Cordovil and Fachada have shown
that the braid monodromy of the complexification
is determined by an associated (unbraided) wiring diagram, and have defined
the braid monodromy of an abstract wiring diagram.
Hironaka's technique may also be applied in this generality.
The algorithm presented here generalizes both these methods.

\subhead 1.2
\endsubhead
Before specializing to arrangements, we present a new interpretation
of the process by which the braid monodromy of a curve $\CC$ is
defined.  This follows in spirit the approach in
\cite{L1}, but uses a self-contained argument based on Hansen's theory
of polynomial covering maps, \cite{H1}, \cite{H2}.  Given a simple
Weierstrass polynomial $f:X\times \C\to \C$ of degree $n$, we consider
the space $Y=X\times \C \setminus \{f(x,z)=0\}$.
In Theorem~2.3, we show that the projection
$p=\pr_1|_Y: Y\to X$ is a fiber bundle map, with structure group
the braid group $B_n$, and monodromy the homomorphism from
$\pi_1(X)$ to $B_n$ induced by the coefficient map of $f$.

This result is applied in the situation where
$f$ defines a plane curve $\CC$, and
$X=\C\setminus \{y_1,\dots, y_s\}$ is the set of regular values of a
generic linear projection.  The braid monodromy of $\CC$ is simply the
coefficient homomorphism, $\a:F_s\to B_n$.   This map depends on
choices of projection, generating curves, and basepoints.
However, the braid-equivalence class of the monodromy---the double
coset $[\a]\in B_s\backslash \Hom(F_s,B_n)\slash B_n$,
where $B_s$ acts on the left by the Artin representation, and $B_n$ acts on
the right by conjugation---is uniquely determined by $\CC$.

For a line arrangement $\A$, changes in the various choices noted
above give rise to changes in the associated braided wiring diagram $\W$.
These, and other, ``Markov moves'' do not affect the braid monodromy.
In practice, the braided wiring diagram of a given arrangement may be
simplified via these moves.  Such simplifications, together with use
of the braid relations, make the braid monodromy presentation of
the group of a complex arrangement accessible.  Furthermore, braided
wiring diagrams associated to arrangements which are lattice-isotopic
in the sense of Randell \cite{R2} are related by Markov moves.  A
combinatorial characterization
of this fact remains to be determined.  Such a characterization,
suggested for (unbraided) wiring diagrams by
Bj\" orner, Las Vergnas, Sturmfels, White, and Ziegler
in \cite{BLSWZ}, Exercise 6.12,
would likely lead to the development of a Jones-type polynomial
for arrangements.

The braid monodromy is also useful in defining Alexander-type
invariants of plane algebraic curves.  Given a curve $\CC$ with
braid monodromy $\a:F_s\to B_n$, one may consider a representation
$\theta: B_n\to \GL(N,R)$, and compute the module of coinvariants
of $\theta\circ \a$.  As noted by Libgober in \cite{L3}, the
$R$-module $A_{\theta}(\CC)=H_0(F_s; R^N_{\theta\circ \a})$
depends only on the equisingular isotopy class of $\CC$ (and on $\theta$).
When $\theta$ is the Burau representation, $A_{\theta}(\CC)$
equals the Alexander module, and thus depends only on
$\pi_1(\C^2\setminus \CC)$.  For other representations
of the braid group, such as the generalized Burau representations
of \cite{CS1}, the module $A_{\theta}(\CC)$ is more likely to be a
homeomorphism-type (rather than homotopy-type) invariant of
the complement, see the discussion in~1.3, and section~7.
For a detailed analysis of Alexander invariants of hyperplane
arrangements, based on the techniques developed in this paper,
we refer to \cite{CS2}.

\subhead 1.3
\endsubhead
In general, the braid monodromy of a plane algebraic curve depends
not only on the number and type of singularities,
but on the relative positions of the singularities as well.
A famous example of Zariski \cite{Z1}, \cite{Z2} consists
of two sextics, both with six cusps, one with all cusps on a conic,
the other not.  Explicit braid monodromy generators for these curves
were given by Rudolph \cite{Ru}, Example 3.  As shown by Zariski, the
two curves have distinct fundamental groups.  Further information
concerning such ``Zariski couples'' may be found in \cite{A-B}.
An example of a different nature is given in~7.4.  There, the two
sextics have the same number of double points (9) and triple
points (2); their fundamental groups are isomorphic, but,
nevertheless, their braid monodromies are not braid-equivalent.

The above example provides an affirmative answer to a question
of Libgober, who raised the possibility in \cite{L3} that the braid
monodromy of a plane algebraic curve which is transverse to the line
at infinity carries more information than the
fundamental group of the complement.  The sextics in~7.4 define
arrangements, originally studied by Falk \cite{Fa}, with distinct
lattices.  This explains the difference in the braid monodromies:  In
Theorem~7.2, we show that the braid-equivalence class of the
monodromy of an arrangement determines the lattice.  On the other
hand, as Falk demonstrated with these examples, the homotopy type of
the complement of an arrangement does not determine the lattice.
However, as noted by Jiang and Yau \cite{JY}, the
complements of these arrangements are not homeomorphic.  This,
and other evidence, suggests that the braid monodromy of a
curve is more closely tied to the homeomorphism type of the
complement (or even to the ambient homeomorphism type of the curve)
than to the fundamental group of the complement.

In the other direction, using classical configurations of MacLane
\cite{MacL}, Rybnikov \cite{Ry} constructs complex arrangements with
isomorphic lattices and distinct fundamental groups.  It follows
that the lattice of a complex arrangement does not determine the
braid monodromy.  We provide another illustration of this phenomenon.
In Theorem~3.9, we show that complex conjugate algebraic
curves have {\it equivalent} braid monodromies.
However, we show in~7.7 that the monodromies of a pair of
conjugate arrangements associated to MacLane's configurations
are not {\it braid-equivalent}, despite the fact that these arrangements
have isomorphic lattices and groups (and, in fact, diffeomorphic complements).

It is not known whether the lattice of a real arrangement determines the
braid monodromy of its complexification.  A result along these lines
may be found in \cite{CF}.  There, the image in the pure braid group
of the braid monodromy of a wiring diagram $\W$ is called the braid monodromy
group of $\W$.  Cordovil and Fachada show that wiring diagrams which
determine identical matroids give rise to equal braid monodromy groups.
This result is not as widely applicable as it may appear.  In~7.5,
we consider arrangements with isomorphic (oriented) matroids and
homeomorphic complements.  Their monodromies are braid-equivalent,
but the associated braid monodromy groups are not conjugate
subgroups of the pure braid group.

\remark{Conventions}  Given elements $x$ and $y$ in a group $G$,
we will write $x^y=y^{-1}xy$ and $[x,y]=xyx^{-1}y^{-1}$.  Also,
we will denote by $\Aut(G)$ the group of {\it right} automorphisms
of $G$, with multiplication $\a\cdot\b=\beta\circ \a$.
\endremark

\remark{Acknowledgments}
We would like to thank Mike Falk for useful conversations, and for
sharing with us his unpublished work with Bernd Sturmfels, as well
as that of John Keaty.
\endremark

\head 2. Polynomial covers and $B_n$-bundles
\endhead

We begin by reviewing polynomial covering maps.  These
were introduced by Hansen in \cite{H1}, and studied in detail in
his book \cite{H2}, which, together with 
Birman's book
\cite{Bi}, is our basic
reference for this section.  We then consider bundles whose
structure group is Artin's braid group $B_n$, and relate them
to polynomial $n$-fold covers.

\subhead 2.1 Polynomial covers
\endsubhead
Let $X$ be a path-connected space that has the homotopy
type of a CW-complex.  A {\it simple Weierstrass polynomial} of
degree $n$ is a map $f: X\times \C \to \C$ given by
$$f(x,z)=z^n+\sum_{i=1}^n a_i(x) z^{n-i},
$$
with continuous coefficient maps $a_i:X\to \C$, and
with no multiple roots for any $x \in X$.
Given such $f$, the restriction of the
first-coordinate projection map $\pr_1:X \times \C \to X$ to the
subspace
$$E=E(f)=\{ (x,z)\in X\times \C \mid f(x,z)=0\}
$$
defines an $n$-fold cover $\pi=\pi_f:E \to X$, the
{\it polynomial covering map} associated to $f$.

Since $f$ has no multiple roots, the {\it coefficient map}
$a=(a_1,\dots ,a_n):X \to \C^n$ takes values in the complement
of the discriminant set, $B^n=\C^n \setminus \Delta_n$.
Over $B^n$, there is a canonical $n$-fold polynomial
covering map $\pi_n=\pi_{_{f_n}}: E(f_n) \to B^n$, determined by the
Weierstrass polynomial $f_n(x,z)=z^n + \sum_{i=1}^n x_i z^{n-i}.$
Clearly, the polynomial cover $\pi_f:E(f)\to X$ is the
pull-back of $\pi_n:E(f_n) \to B^n$ along the coefficient map
$a: X\to B^n$.

This can be interpreted on the level of fundamental groups
as follows.  The fundamental group of the configuration space,
$B^n$, of $n$ unordered points in $\C$ is the group, $B_n$,
of braids on $n$ strands.
The map $a$ determines the {\it coefficient homomorphism}
$\a=a_*: \pi_1(X) \to B_n$, unique up to conjugacy.  One
may characterize polynomial covers as those covers $\pi: E\to X$
for which the characteristic homomorphism to the symmetric group,
$\chi:\pi_1(X)\to \Sigma_n$, factors through the canonical surjection
$\tau_n: B_n \to \Sigma_n$ as $\chi = \tau_n \circ \a.$

Now assume that the simple Weierstrass polynomial $f$
is {\it completely solvable}, that is, factors as
$$f(x,z)=\prod_{i=1}^n (z-b_i(x)),
$$
with continuous roots $b_i: X \to \C$.  Since the Weierstrass
polynomial $f$ is simple, the {\it root map} $b=(b_1,\dots ,b_n):X\to
\C^n$ takes values in the complement, $P^n=\C^n \setminus \A_n$, of the
braid arrangement
$\A_n=\{\ker (w_i-w_j)\}_{1\le i<j\le n}$.
Over $P^n$, there is a canonical $n$-fold covering map,
$q_n=\pi_{_{Q_n}}: E(Q_n) \to P^n$, determined by the Weierstrass
polynomial $Q_n(w,z)=(z-w_1)\cdots (z-w_n)$.
Evidently, the cover $\pi_f:E\to X$ is the pull-back of
$q_n: E(Q_n)\to P^n$ along the root map $b: X\to P^n$.

The fundamental group of the configuration space,
$P^n$, of $n$ ordered points in $\C$ is the group, $P_n=\ker\tau_n$,
of pure braids on $n$ strands.
The map $b$ determines
the {\it root homomorphism} $\b=b_*: \pi_1(X) \to P_n$,
unique up to conjugacy.
The polynomial covers which are trivial covers (in the usual sense)
are precisely those for which the coefficient homomorphism
factors as $\a=\iota_n\circ \b$,
where $\iota_n: P_n\to B_n$ is the canonical injection.

\subhead 2.2   $B_n$-Bundles
\endsubhead
The group $B_n$ may be realized as the mapping class group 
$\frak M_{0,1}^n$ of orientation-preserving diffeomorphisms of 
the disk $D^2$, permuting a collection of $n$ marked points.  
Upon identifying $\pi_1(D^2 \setminus \{ n \text{ points}\})$ with 
the free group $F_n$, the action of $B_n$ on $\pi_1$ yields the 
{\it Artin representation}, $\a_n:B_n\to \Aut(F_n)$.  As showed by 
Artin, this representation is faithful.  Hence, we 
may---and often will---identify a braid $\theta\in B_n$ with the 
corresponding braid automorphism, $\a_n(\theta)\in \Aut(F_n)$.

Now let $f:X\times \C\to \C$ be a simple Weierstrass polynomial.
Let $\pi_f: E(f)\to X$ be the corresponding polynomial
$n$-fold cover, and $a:X \to B^n$ the coefficient map.
Consider the complement
$$Y=Y(f)=X\times \C \setminus E(f),$$
and let $p=p_f: Y(f)\to X$ be the restriction of
$\pr_1:X\times \C \to X$ to $Y$.

\proclaim{Theorem 2.3}
The map $p:Y\to X$ is a locally trivial bundle, with structure group
$B_n$ and fiber $\C_n=\C \setminus \{ n \text{ points}\}$.  Upon
identifying $\pi_1(\C_n)$ with $F_n$, the monodromy
of this bundle may be written as $\a_n\circ \a$, where
$\a=a_*:\pi_1(X)\to B_n$ is the coefficient homomorphism.

Moreover, if $f$ is completely solvable, the structure group
reduces to $P_n$, and the monodromy factors as $\a_n\circ\iota_n\circ \b$,
where $\b=b_*:\pi_1(X)\to P_n$ is the root homomorphism.
\endproclaim

\demo{Proof}  We first prove the theorem for the configuration
spaces, and their canonical Weierstrass polynomials.
Start with $X=P^n$, $f=Q_n$, and the canonical cover $q_n: E(Q_n)\to P^n$.
Clearly, $Y(Q_n)=\C^{n+1} \setminus E(Q_n)$ is equal to the configuration
space $P^{n+1}$.
Let $\rho_n=p_{_{Q_n}}:P^{n+1}\to P^n$ be the restriction of
$\pr_1:\C^n\times \C\to \C^n$.  As shown by Faddell and Neuwirth~\cite{FN},
this is a bundle map, with fiber $\C_n$, and monodromy the restriction of
the Artin representation to $P_n$.

Next, consider $X=B^n$, $f=f_n$, and the canonical cover
$\pi_n:E(f_n) \to B^n$. Forgetting the order of the points defines
a covering projection from the ordered to the unordered configuration
space, $\kappa_n:P^n \to B^n$.
In coordinates, $\kappa_n(w_1,\dots ,w_n)=(x_1,\dots ,x_n)$, where
$x_i=(-1)^i s_i(w_1,\dots ,w_n)$, and $s_i$ are the elementary
symmetric functions.  By Vieta's formulas, we have
$$Q_n(w,z) = f_n(\kappa_n(w),z).$$
Let $Y^{n+1}=Y(f_n)$ and $p_n=p_{f_n}: Y^{n+1} \to B^n$.
By the above formula, we see that
$\kappa_n\times \id :P^n\times \C \to B^n\times \C$ restricts
to a map $\bar\kappa_{n+1}:Y(Q_n) \to Y(f_n)$, which fits into
the fiber product diagram
$$
\CD
P^{n+1}         @>\rho_n>>      P^n \\
@VV\bar\kappa_{n+1}V      @VV{\kappa_n}V \\
Y^{n+1}         @>p_n>>      B^n
\endCD
$$
where the vertical maps are principal $\Sigma_n$-bundles.
Since the bundle map $\rho_n: P^{n+1}\to P^n$ is equivariant with
respect to the $\Sigma_n$-actions, the map on quotients,
$p_n: Y^{n+1} \to B^n$, is also a bundle map, with fiber $\C_n$,
and monodromy action the Artin representation of $B_n$.
This finishes the proof in the case of the canonical
Weierstrass polynomials over configuration spaces.

Now let $f:X\times \C\to \C$ be an arbitrary simple Weierstrass
polynomial.  We then have the following cartesian square:
$$
\CD
Y      @>>>   Y^{n+1} \\
@VVpV           @VVp_nV \\
X      @>a>>   B^n
\endCD
$$
In other words, $p: Y\to X$ is the pullback of the bundle
$p_n: Y^{n+1}\to B^n$ along the coefficient map $a$.
Thus, $p$ is a bundle map, with fiber $\C_n$, and monodromy
representation $\a_n\circ \a$.  When $f$ is completely
solvable, the bundle $p: Y\to X$ is the pullback of
$\rho_n: P^{n+1}\to P^n$ along the root map $b$.
Since $\a=\iota_n\circ \b$, the monodromy is as claimed.
\qquad\qed
\enddemo

\example{Remark 2.4}
Let us summarize the above discussion of braid bundles over
configuration spaces.  From the Faddell-Neuwirth
theorem, it follows that $P^n$ is a $K(P_n,1)$ space.  Since
the pure braid group is discrete, the classifying $P_n$-bundle
(in the sense of Steenrod) is the universal cover
$\widetilde{P}^n \to P^n$.  We considered two bundles
over $P^n$, both associated to this one:
\roster
\item"(i)" $q_n:E(Q_n) \to X^n$, by the trivial representation of
$P_n$ on $\{1,\dots ,n\}$;
\item"(ii)" $\rho_n:P^{n+1} \to P^n$, by the (geometric) Artin
representation of $P_n$ on $\C_n$.
\endroster
Since $B^n$ is covered by $P^n$, it is a $K(B_n,1)$ space, and the
classifying $B_n$-bundle is $\widetilde{B}^n \to B^n$.  There
were three bundles over $B^n$ that we mentioned, all associated
to this one:
\roster
\item"(iii)" $\kappa_n:X^n \to B^n$, by the canonical surjection
$\tau_n: B_n\to \Sigma_n$;
\item"(iv)" $\pi_n:E(f_n) \to B^n$, by the above, followed by the
permutation representation of $\Sigma_n$ on $\{1,\dots ,n\}$;
\item"(v)" $p_n:Y^{n+1} \to B^n$, by the (geometric) Artin
representation of $B_n$ on $\C_n$.
\endroster
Finally, note that $\pi_1(Y^{n+1})$ is isomorphic to
$B_n^1=F_n\rtimes_{\a_n} B_n$, the group of braids on $n+1$
strands that fix the endpoint of the last strand, and that
$Y^{n+1}$ is a $K(B_n^1,1)$ space.
\endexample

\head 3.  The braid monodromy of a plane algebraic curve
\endhead

We are now ready to define the braid monodromy of an algebraic
curve in the complex plane.  The construction, based on
classical work of Zariski and van~Kampen, is due to
Moishezon \cite{Mo}.
We follow the exposition of Libgober \cite{L1}, \cite{L2},
\cite{L3}, but interpret the construction in the context
established in the previous section.

\subhead 3.1 The construction
\endsubhead
Let $\CC$ be a reduced algebraic curve in $\C^2$, with defining polynomial
$f$ of degree $n$.  Let $\pi:\C^2\to \C$ be a linear projection,
and let $\Y =\{y_1,\dots, y_s\}$ be the set of points in $\C$ for which
the fibers of $\pi$ contain singular points of $\CC$, or are tangent
to $\CC$.  Assume that $\pi$ is
generic with respect to $\CC$.  That is, for each $k$, the line
$\L_k=\pi^{-1}(y_k)$ contains at most one singular point $v_k$ of $\CC$, and
does not belong to the tangent cone of $\CC$ at $v_k$,
and, moreover, all tangencies are simple.  Let $\L$ denote the union
of the lines $\L_k$, and let $y_0$ be a basepoint in $\C\setminus \Y$.
The definition of the braid monodromy of $\CC$ depends on two observations:

\proclaim\nofrills{(i)}\ \ The restriction of the projection map,
$p:\C^2 \setminus \CC \cup \L \to \C \setminus \Y$,
is a locally trivial bundle.
\endproclaim

Fix the fiber $\C_n=p^{-1}(y_0)$ and a basepoint $\hat y_0 \in \C_n$.
The monodromy of $\CC$ is, by definition, the holonomy of this bundle,
$\rho: \pi_1(\C \setminus \Y,y_0) \to \Aut(\pi_1(\C_n,\hat y_0))$.
Upon identifying $\pi_1(\C \setminus \Y,y_0)$ with $F_s$,
and $\pi_1(\C_n,\hat y_0)$ with $F_n$, this can be written
as $\rho: F_s \to \Aut(F_n)$.

\proclaim\nofrills{(ii)}\ \ The image of $\rho$ is contained in the
braid group $B_n$ (viewed as a subgroup of $\Aut(F_n)$ via the Artin
embedding $\a_n$).
\endproclaim

The {\it braid monodromy} of $\CC$ is the homomorphism $\a:F_s\to B_n$
determined by $\a_n\circ \a=\rho$.

We shall present a self-contained proof of these two assertions,
and, in the process, identify the map $\a$.  The first
assertion is well-known, and can also be proved by standard techniques
(using blow-ups and Ehresmann's criterion---see \cite{Di}, page~123),
but we find our approach sheds some light on the underlying
topology of the situation.

\subhead 3.2  Braid monodromy and polynomial covers
\endsubhead
Let $\pi:\C^2 \to \C^1$ be
a linear projection, generic with respect to the given
algebraic curve $\CC$ of degree $n$.  We may assume
(after a linear change of variables in $\C^2$ if necessary) that
$\pi=\pr_1$, the projection map onto first coordinate.
In the chosen coordinates, the defining polynomial $f$ of $\CC$
may be written as
$f(x,z)=z^n+\sum_{i=1}^n a_i(x) z^{n-i}$.
Since $\CC$ is reduced, for each $x\notin \Y$, the equation $f(x,z)=0$
has $n$ distinct roots.  Thus, $f$ is a simple Weierstrass polynomial
over $\C\setminus \Y$, and
$$\pi=\pi_{f}:\CC \setminus \CC\cap \L \to \C\setminus \Y\tag{1}$$
is the associated polynomial $n$-fold cover.  Note that
$Y(f)=((\C \setminus \Y)\times \C) \setminus (\CC \setminus
\CC\cap \L)=\C^2 \setminus \CC \cup \L$.  By Theorem~2.3,
the restriction of $\pr_1$ to $Y(f)$,
$$p:\C^2 \setminus \CC \cup \L \to \C \setminus \Y,\tag{2}$$
is a bundle map, with structure group $B_n$, fiber $\C_n$, and
monodromy homomorphism
$$\a=a_*:\pi_1(\C\setminus \Y) \to B_n. \tag{3}$$
This proves assertions (i) and (ii).  Furthermore, we have

\proclaim{Theorem 3.3} The braid monodromy of a plane algebraic
curve coincides with the coefficient homomorphism of the
associated polynomial cover.
\endproclaim

In the case where $\CC=\A$ is an arrangement of (affine) lines
in $\C^2$, more can be said.  First, the critical set
$\Y=\{y_1,\dots ,y_s\}$ consists (only) of the images under $\pi=\pr_1$
of the vertices of the arrangement.  Furthermore, a defining
polynomial for $\A$ can be written as
$f(x,z)=\prod_{i=1}^n (z-\ell_i(x))$,
where each $\ell_i$ is a linear function in $x$.
Thus, the associated polynomial cover is trivial, and the
monodromy representation is
$$\lambda=\ell_*: \pi_1(\C\setminus \Y) \to P_n.
$$
An explicit formula for $\lambda$ will be given in section~5.
For now, let us record the following:

\proclaim{Theorem 3.4} The pure braid monodromy of a line arrangement
coincides with the root homomorphism of the associated (trivial)
polynomial cover.
\endproclaim

\subhead 3.5  Braid equivalence
\endsubhead
The braid monodromy of a plane algebraic curve is
not unique, but rather, depends on the choices made in defining it.
This indeterminacy was studied by Libgober in \cite{L2}, \cite{L3}.
To make the analysis more precise, we first need a definition.

\definition{Definition 3.6}  Two homomorphisms $\a:F_s\to B_n$ and
$\a':F_s\to B_n$ are {\it equivalent} if there exist
automorphisms $\psi\in \Aut(F_s)$ and $\phi\in \Aut(F_n)$ with
$\phi(B_n)\subset B_n$ such that
$\a'(\psi(g)) = \phi^{-1}\cdot \a(g)\cdot \phi$,
for all $g\in F_s$.  In other words, the following diagram commutes
$$
\CD
F_s        @>\a>>      B_n \\
@VV{\psi}V      @VV\conj_{\phi}V \\
F_s        @>\a'>>     B_n
\endCD
$$
If, moreover, $\psi\in B_s$ and $\phi\in B_n$, the homomorphisms
$\a$ and $\a'$ are {\it braid-equivalent}.
\enddefinition

\proclaim{Theorem 3.7}  The braid monodromy of a plane algebraic
curve $\CC$ is well-defined up to braid-equivalence.
\endproclaim

\demo{Proof}
First fix the generic projection.
The identification $\pi_1(\C\setminus \Y)=F_s$ depends on the
choice of a ``well-ordered'' system of generators (see \cite{Mo}
or the discussion in 4.1), and any two such choices
yield monodromies which differ by a braid automorphism of $F_s$,
see \cite{L2}.  Furthermore, there is the choice of basepoints,
and any two such choices yield monodromies differing by a conjugation
in $B_n$.

Finally, one must analyze the effect of a change in the choice of
generic projection.   Let $\pi$ and $\pi'$ be two such projections,
with critical sets $\Y$ and $\Y'$, and braid monodromies
$\a:\pi_1(\C\setminus \Y)\to B_n$ and $\a':\pi_1(\C\setminus \Y')\to B_n$.
Libgober \cite{L3} shows that there is a homeomorphism
$h:\C\to\C$, isotopic to the identity, and taking $\Y$ to $\Y'$,
for which the isomorphism
$h_*:\pi_1(\C \setminus \Y) \to \pi_1(\C \setminus \Y')$
induced by the restriction of $h$ satisfies $\a'\circ h_* =\a$.
From the construction, we see that $h$ can be taken to be the
identity outside a ball of large radius (containing $\Y\cup \Y'$).
Thus, once the identifications of source and target with $F_s$ are made,
$h_*$ can be written as the composite of an inner automorphism of $F_s$
with a braid automorphism of $F_s$:  $h_*=\conj_{g}\circ \psi$.
Trading the inner automorphism of $F_s$ for an inner
automorphism of $B_n$, we obtain $\a'\circ \psi=\conj_{\a'(g)}\circ \a$,
completing the proof.\qquad \qed
\enddemo

Thus, we may regard the braid monodromy of $\CC$ as a braid-equivalence
class, i.e., as a double coset
$[\a]\in B_s\backslash \Hom(F_s,B_n)\slash B_n$, uniquely determined
by $\CC$.  In fact, it follows from \cite{L3} that $[\a]$ depends only
on the equisingular isotopy class of the curve.

\subhead 3.8  Conjugate curves
\endsubhead
If $\CC$ is a plane curve with defining polynomial $f=f(x,z)$
of degree $n$, let $\bar\CC$ be the curve defined by the polynomial
$\bar f$  whose coefficients are the complex conjugates of those of
$f$. In other words, $\bar f(x,z)= \overline{f(\bar x,\bar z)}$.
In this subsection, we relate the braid monodromies of $\CC$ and
$\bar\CC$.  In general, the braid monodromies of conjugate curves are
not braid-equivalent, as shown in 7.7.  Nevertheless, we have the
following:

\proclaim{Theorem~3.9}  The braid monodromies of conjugate curves are
equivalent.
\endproclaim

\demo{Proof}  Let $\CC$ and $\bar\CC$ be conjugate curves defined by
polynomials $f$ and $\bar f$ of degree
$n$.  Choose coordinates in $\C^2$ so that $\pi=\pr_1$
is generic with respect to $\CC$.
Then $\pi$ is evidently also generic with respect to $\bar\CC$.
Let $\Y$ and $\bar \Y$ be the critical sets of $\CC$ and $\bar\CC$
with respect to this projection.  Complex conjugation $\C\to\C$
restricts to a map $d:\C\setminus\Y \to \C\setminus\bar \Y$.
Choose a basepoint $y_0$ with $\Imag(y_0)=0$.  Then $d$ induces
an isomorphism $d_*:\pi_1(\C\setminus\Y,y_0) \to
\pi_1(\C\setminus\bar \Y,y_0)$.  Identifying these groups with
$F_s=\langle x_1,\dots,x_s\rangle$, we have $d_*=\delta_s$, where
$\delta_s \in \Aut(F_s)$ is given by
$\delta_s(x_k) = (x_1\cdots x_{k-1})
\cdot x_k^{-1} \cdot (x_1\cdots x_{k-1})^{-1}.$

Since the discriminant locus $\Delta_n$ is defined by
a polynomial with real coefficients, complex conjugation
$\C^n \to \C^n$ restricts to a map $e:B^n \to B^n$.
The induced map $\epsilon_n = e_*:B_n\to B_n$ is readily
seen to be the automorphism defined on generators by
$\epsilon_n(\sigma_i)=\sigma_i^{-1}$.  As shown by Dyer and Grossman
\cite{DG}, this involution generates $\Out(B_n)=\Z_2$, for $n\ge 3$.

Let $a$ and $\bar a$ be the coefficient maps of $f$ and $\bar f$ respectively.
The fact that the defining polynomials of $\CC$ and $\bar\CC$ have
complex conjugate coefficients may be expressed as
$\bar a \circ d = e\circ a$.
Passing to fundamental groups, we have
$\bar \a \circ \delta_s = \epsilon_n\circ \a$.  Checking that
$\epsilon_n =\conj_{\delta_n}$ (see \cite{DG})
completes the proof.\quad\qed
\enddemo

\head 4.  The fundamental group of a plane algebraic curve
\endhead
We now give the braid monodromy presentation
of the fundamental group of the complement of a plane algebraic
curve $\CC$. This presentation first appeared in the classical work
of van~Kampen and Zariski \cite{vK}, \cite{Z2}, and has been much
studied since, see e.g. \cite{Mo}, \cite{MT}, \cite{L1}, \cite{L2},
\cite{Ru}, \cite{Di}.

\subhead 4.1 Braid monodromy presentation
\endsubhead
The homotopy exact sequence of the bundle
$p:\C^2 \setminus \CC \cup \L \to \C \setminus \Y$
of \thetag{2} reduces to
$$1\to \pi_1(\C_n) \to \pi_1(\C^2 \setminus \CC \cup \L)
@>p_*>>\pi_1(\C \setminus \Y)\to 1.
$$
This sequence is split exact, with action given by
the braid monodromy homomorphism $\a$ of \thetag{3}.
To extract a presentation of the middle group,
order the points of $\Y$ by decreasing real part,
and pick the basepoint $y_0$ in $\C \setminus \Y$ with
$\Real(y_0)>\max \{\Real(y_k)\}$.
Choose loops $\xi_k:[0,1]\to\C \setminus \Y$ based at $y_0$,
going up and above $y_1,\dots, y_{k-1}$, passing around $y_k$
in counterclockwise direction, and coming back the same way.
Setting $x_k=[\xi_k]$, identify $\pi_1(\C \setminus \Y, y_0)$
with $F_s=\langle x_1,\dots ,x_s\rangle$.  Similarly, identify
$\pi_1(\C_n, \hat y_0)$ with $F_n=\langle t_1,\dots , t_n\rangle$.
Having done this, $\pi_1(\C^2 \setminus \CC \cup \L, \hat y_0)$
becomes identified with the semidirect product
$F_n \rtimes_{\a} F_s$. The corresponding presentation is
$$\pi_1(\C^2 \setminus \CC \cup \L)=
\langle t_1,\dots t_n,x_1\dots ,x_s \mid
x_k^{-1}\cdot t_i\cdot x_k = \a(x_k)(t_i)\rangle.
$$

The fundamental group of the complement of the curve is the
quotient of $\pi_1(\C^2 \setminus \CC \cup \L)$
by the normal closure of $F_s=\langle x_1,\dots ,x_s\rangle$.
Thus, $\pi_1(\C^2 \setminus \CC)=
\langle t_1,\dots ,t_n\mid t_i = \a(x_k)(t_i)\rangle.$
This presentation can be simplified by Tietze-II
moves---eliminating redundant relations.  Doing so,
one obtains the
{\it braid monodromy presentation}
$$\pi_1(\C^2 \setminus \CC)=
\langle t_1,\dots ,t_n\mid t_i = \a(x_k)(t_i),
\ \ i=j_1,\dots , j_{m_k-1} ; \ \ k=1,\dots, s
\rangle.\tag{4}$$
If $y_k$ corresponds to a singular point of $\CC$, then $m_k$
denotes the multiplicity of that singular point, while if $y_k$
corresponds to a (simple) tangency point, $m_k=2$.  In either
case, the indices $j_1,\dots,j_{m_k-1}$ must be chosen  
appropriately, see \cite{L1} and the discussions in 5.1 and 6.1.

Let $K(\CC)$ be the $2$-complex modeled on the braid monodromy
presentation.  There is an obvious embedding of
this complex into the complement of $\CC$.  The main result
of \cite{L1} is the following.

\proclaim{Theorem 4.2 (Libgober)}
The $2$-complex $K(\CC)$ is a deformation retract of
$\C^2 \setminus \CC$.
\endproclaim

\example{Remark 4.3}  The group $G(\a)$ defined by presentation
\thetag{4} is the quotient of $F_n$ by the normal subgroup generated
by $\{\gamma(t)\cdot t^{-1} \mid \gamma\in \im(\a),\ t\in F_n\}$.
In other words,
$G(\a)$ is the maximal quotient of $F_n$
on which the representation $\a:F_s \to B_n$ acts trivially.
If $\a':F_s \to B_n$ is equivalent to $\a$, then $G(\a)$ is
isomorphic to $G(\a')$.  Indeed, the equivalence condition
$\a'\circ \psi = \conj_{\phi}\circ \a$ can be written as
$\phi(\a(g)(t)\cdot t^{-1})=\a'(\psi(g))(\phi(t))\cdot
\phi(t)^{-1}$, $\forall g\in F_s,\ \forall t\in F_n$.
Thus $\phi\in \Aut(F_n)$ induces an isomorphism
$\bar\phi:G(\a)\to G(\a')$.
\endexample

\subhead 4.4  Braid monodromy generators
\endsubhead
We now make the presentation \thetag{4} more precise.
First recall that the braid group $B_n$ has generators
$\sigma_1,\dots ,\sigma_{n-1}$, and relations
$\sigma_i\sigma_{i+1}\sigma_i=\sigma_{i+1}\sigma_i\sigma_{i+1}
\ (1\le i<n-1),\ \sigma_i\sigma_j=\sigma_j\sigma_i
\ (|i-j|>1)$, see \cite{Bi}, \cite{H2}.  The Artin representation
$\a_n:B_n\to \Aut(F_n)$ is given by:
$$
\sigma_i(t_j)=\cases  t_it_{i+1}t_i^{-1}&\text{ if } j=i,\\
t_i&\text{ if } j=i+1,\\
t_j &\text{ otherwise.}
\endcases
$$

For each $k=1,\dots, s$, let $\gamma_k\in B_{m_k}<B_n$ be the
``local monodromy'' around $y_k$.  Then
$$\a(x_k)=\beta_k^{-1}\gamma_k \beta_k,$$
where $\beta_k\in B_n$ is the monodromy along the portion of $\xi_k$
from $y_0$ to just before $y_k$.  One would like to express these
braids in terms of the standard generators $\sigma_i$ of $B_n$.
This may be accomplished in two steps.

\remark{Step 1}
The structure of the (isolated) singularity $v_k$ above $y_k$ determines
the local braid $\gamma_k$.  This braid may be obtained from the
Puiseux series expansion of the defining polynomial $f(x,z)$ of $\CC$.
This is implicit in the work of Brieskorn and Kn\"orrer \cite{BK} and
Eisenbud and Neumann \cite{EN}.
\endremark

\example{Example 4.5}  Consider the plane curve $\CC: z^p - x^q = 0$.
The fundamental group of its complement was determined by Oka \cite{Ok}.
A look at Oka's computation reveals that the braid monodromy generator
is $(\sigma_1\cdots \sigma_{p-1})^q\in B_p.$
For instance, to a simple tangency corresponds
$\sigma_1$, to a node, $\sigma_1^2$, and to a cusp, $\sigma_1^3$.
\endexample

\example{Example 4.6} By the above, the braid monodromy generator of a central
line arrangement $\A: z^n-x^n=0$ is a full twist on $n$ strands,
$\Delta^2=(\sigma_1\cdots \sigma_{n-1})^n\in B_n $ (see also \cite{Hir}).
\endexample

\remark{Step 2}
The conjugating braids $\b_k$ depend on the relative positions of
the singularities of $\CC$.  These braids may be specified as follows.
Let $\eta_k$ denote the portion of the path $\xi_k$ from $y_0$ to
just before $y_k$.  The braid $\b_k$ is identified by tracking the
components of the fiber of the polynomial cover
$\pi=\pi_f:\CC\setminus \CC\cap\L \to \C\setminus\Y$ of \thetag{1}
over $\eta_k$.
Generically, these components have distinct real parts.
Braiding occurs when the real parts of two components coincide.
We record this braiding by analyzing
the imaginary parts of the components,
as indicated in the figure below.
\endremark

\bigskip

\centerline{
\beginsmallgraph
\vertex at (2.5, 15)
\vertex at (2, 13.75)
\vertex at (3.25, 16)
\vertex at (23.5, 16)
\vertex at (22, 13.75)
\vertex at (21.5, 13)
\label {$1$} at (21.5, 12.5)
\label {$2$} at (22.5, 13.75)
\label {$3$} at (23.5, 15.5)
\edge from (2.5, 15) to (23.5, 16)
\edge from (2, 13.75) to (15.75, 13.75)
\edge from (18.25, 13.75) to (22, 13.75)
\edge from (3.25, 16) to (6.5, 15.47)
\edge from (8.75, 15.1) to (21.5, 13)
\edge from (21, 11) to (24, 13)
\edge from (21, 11) to (21, 15)
\edge from (21, 15) to (24, 17)
\edge from (24, 13) to (24, 17)
\edge from (1, 11) to (4, 13)
\edge from (1, 11) to (1, 15)
\edge from (1, 15) to (4, 17)
\edge from (4, 13) to (4, 17)
\label {Re$(z)$} at (25, 8.5)
\label {Im$(z)$} at (22.5, 10.5)
\edge from (1, 8.5) to (4, 8.5)
\edge from (6, 8.5) to (9, 8.5)
\edge from (11, 8.5) to (14, 8.5)
\edge from (16, 8.5) to (19, 8.5)
\edge from (21, 8.5) to (24, 8.5)
\edge from (2.5, 7) to (2.5, 10)
\edge from (7.5, 7) to (7.5, 10)
\edge from (12.5, 7) to (12.5, 10)
\edge from (17.5, 7) to (17.5, 10)
\edge from (22.5, 7) to (22.5, 10)
\vertex at (2, 8)
\label {$2$} at (1.5, 8)
\vertex at (2.5, 9.25)
\label {$3$} at (2, 9.25)
\vertex at (3.75, 9.5)
\label {$1$} at (3.25, 9.5)
\vertex at (7, 8)
\label {$2$} at (6.5, 8)
\vertex at (8.5, 9.75)
\label {$3$} at (8, 9.75)
\vertex at (8.5, 9)
\label {$1$} at (8, 9)
\vertex at (12, 8)
\label {$2$} at (11.5, 8)
\vertex at (13.75, 9.75)
\label {$3$} at (13.25, 9.75)
\vertex at (13.25, 9)
\label {$1$} at (12.75, 9)
\vertex at (17, 8)
\label {$2$} at (16.5, 8)
\vertex at (18.75, 9.75)
\label {$3$} at (18.25, 9.75)
\vertex at (17, 9)
\label {$1$} at (16.5, 9)
\vertex at (22, 8)
\label {$2$} at (21.5, 8)
\vertex at (24, 9.75)
\label {$3$} at (23.5, 9.75)
\vertex at (21, 7.5)
\label {$1$} at (20.5, 7.5)
\label {$\sigma_1^{-1}$} at (17.5, 11.5)
\label {$\sigma_2$} at (7.5, 11.5)
\endgraph
}

\nopagebreak
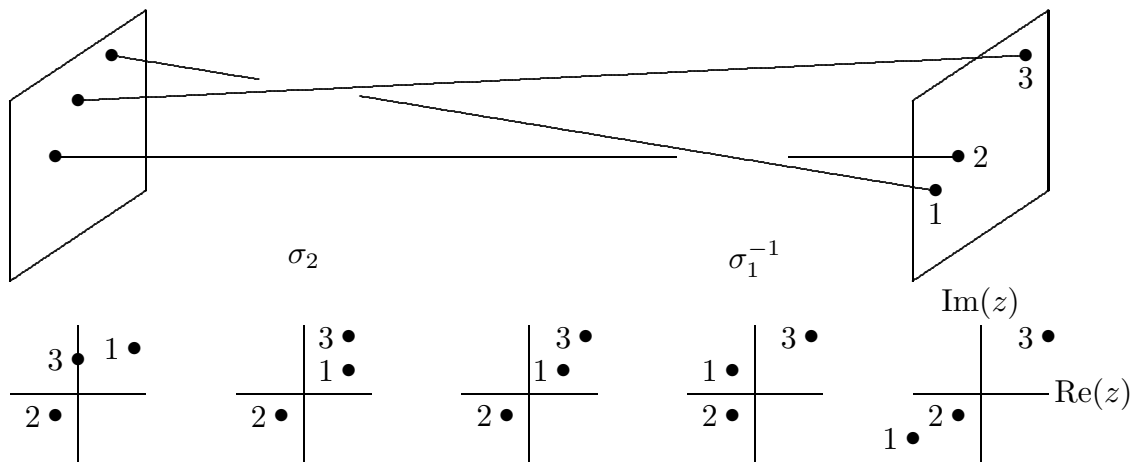
\centerline{Figure 1.  Braiding in a polynomial cover}
\medskip

More explicitly, recall that the polynomial cover $\pi$ is
embedded in the trivial line bundle
$\pr_1:(\C\setminus\Y)\times\C\to \C\setminus\Y$.
Let $y_k'=y_k+\epsilon$ denote the endpoint of the path $\eta_k$.
Without loss of generality, we may assume that the components of
$\pi^{-1}(y_k')$ (resp.~$\pi^{-1}(y_0)$) have distinct real parts.
After an isotopy of $\CC$, we may assume further that the positions of
the components of $\pi^{-1}(y_k')$ in $\pr_1^{-1}(y_k')=\C$ are
identical to those of $\pi^{-1}(y_0)$ in $\pr_1^{-1}(y_0)=\C$.
Then the image of the path $\eta_k$ under the coefficient map
$a:\C\setminus\Y\to B^n$ is a loop $a(\eta_k)$ in configuration space,
and the braid $\b_k$ is the homotopy class of this loop.

\example{Remark 4.7}  The closed braid determined by the product,
$\a(x_1)\cdots \a(x_s)$, of the braid monodromy generators is the link
of the curve $\CC$ at infinity.  In the works of Moishezon and
Libgober,  it is usually assumed that $\CC$ is in general position
relative to the line at infinity in $\CP^2$.  In that situation, the
link at infinity is the $n$-component Hopf link, and thus we have
$\a(x_1)\cdots\a(x_s) =\Delta^2$ by Example~4.6.
\endexample

\head 5. The braid monodromy of a complex arrangement
\endhead

The fundamental group of the complement of an arrangement of complex
hyperplanes is, by a well-known Lefschetz-type theorem of Zariski,
isomorphic to that of a generic two-dimensional section.
So let $\A$ be an arrangement of $n$ complex lines in $\C^2$, with
group $G=\pi_1(\C^2\setminus\A)$.  In this section, we provide
an explicit description of the pure braid monodromy of $\A$.

\subhead 5.1 Braided wiring diagrams
\endsubhead
Choose coordinates in $\C^2$ so that the projection
$\pi=\pr_1:\C^2\to \C$ is generic with respect to $\A$, and let
$f(x,z)=\prod_{i=1}^n (z-\ell_i(x))$ be a defining polynomial for
$\A$.  The points $y_k\in\Y$ are the images under $\pi$ of the
vertices of $\A$.  These vertices, the points
$v_k=(y_k,z_k)\in\C^2$ where $z_k=\ell_{i_1}(y_k)=\dots=
\ell_{i_r}(y_k)$ for $r\ge 2$, are the only singularities of $\A$;
there are no tangencies.  Without loss of generality, assume that the
points $y_k$ have distinct real parts.  As noted in~4.4, {\it Step 1}, the
local monodromy around $y_k$ depends only on $v_k$.  It is completely
determined by the multiplicity of $v_k$, and the relative positions
of the lines incident on $v_k$.  These data, and the braiding of
the lines of $\A$ over the paths $\eta_k$, determine the braid
monodromy of the arrangement.  All of this information may be
effectively recorded as follows.

Order the points of $\Y$ as before, and choose the
basepoint $y_0\in \C\setminus \Y$ so that
$\Real(y_0)>\Real(y_1)>\dots>\Real(y_s)$.
Let $\xi:[0,1]\to\C$ be a (smooth) path emanating from
$y_0$ and passing through $y_1,\dots,y_s$ in order.
Note that we may take the path $\xi$ to be
a horizontal line segment in a neighborhood of each $y_k$.
Call such a path admissible.  Let
$$\W = \{(x,z)\in \xi\times\C\mid f(x,z)=0\}$$
be the {\it braided wiring diagram} associated to $\A$.
Note that $\W$ depends on the generic linear projection
$\pi$ and on the admissible path $\xi$.
If $\{z=\ell_i(x)\}$ is a line of $\A$, we
call $\W \cap \{z=\ell_i(x)\}$ the associated wire.
Since the path $\xi$ passes through the points of $\Y$,
the vertices of $\A$ are contained in $\W$.

Over portions of the path $\xi$ between the points of $\Y$,
the lines of $\A$ (resp.~wires of $\W$) may braid.  Let
$y_k'= y_k+\epsilon$, and $y_k'' = y_k-\epsilon$, for some
sufficiently small $\epsilon$.  We may assume that, over
$y_k'$ and $y_k''$, the wires of $\W$ (i.e.,~the components of
the fiber of the polynomial cover $\pi_f$) have distinct real parts.
Arguing as in~4.4, {\it Step 2}, we associate a braid $\b_{k,k+1}$
to the portion of $\xi$ from $y_k''$ to $y_{k+1}'$.

After an isotopy of $\A$, we may also assume that the
positions of the wires of $\W$ over the points
$y_0$, $y_k'$, and $y_k''$ are all identical.
Thus a braided wiring diagram $\W$ may be abstractly
specified  by a sequence of states, vertices, and braids:
$$\P_{s+1}@<{\ \b_{s,s+1}}<<@<{V_s}<<\P_s@<<<@<<<{\cdots}@<<<@<<<
\P_2@<{\ \b_{1,2}}<<@<{V_1}<<\P_1@<{\ \b_{0,1}}<<\P_0,$$
where the states $\P_k$ are
permutations of $\{1,\dots,n\}$ beginning with the identity
permutation and recording the relative heights of the wires.
The vertex set $V_k=\{i_1,\dots,i_r\}$ records the indices of
the wires incident on the $k^{\text{th}}$ vertex $v_k$ of $\A$ (in
terms of the order given by the initial state $\P_0$).
The braids $\b_{k,k+1}$ are obtained as above.
By choosing the basepoint
$y_0$ sufficiently close to $y_1$, we may assume that the initial
braid $\b_{0,1}$ is trivial.  If such a diagram is depicted
as above, the braids $\b_{k,k+1}$ should be read off from left to right.
Note that the this notion generalizes
that of a wiring diagram due to Goodman \cite{Go},
and that the admissible 2-graphs utilized by Arvola \cite{Ar}, \cite{OT},
may be viewed as examples of braided wiring diagrams.
Explicit examples are given in section~7.

\subhead 5.2  Generators of $P_n$
\endsubhead
Before proceeding, we need to review some facts about the
pure braid group $P_n=\ker(\tau_n:B_n\to \Sigma_n)$.
This group has generators
$$A_{i,j}=\s_{j-1}\cdots \s_{i+1}\cdot \s_i^{2}
\cdot \s_{i+1}^{-1}\cdots \s_{j-1}^{-1},
\quad 1\le i<j\le n,
$$
and relations that set a generator equal to a certain conjugate
of itself, see \cite{Bi}, \cite{Ha}.  In particular,
$H_1(P_n)=\Z^{\binom{n}{2}}$ is generated by the images of the generators
$A_{i,j}$.
The conjugation action of $B_n$ on $P_n$ is given by the
following formulas (compare \cite{DG}):
$$A_{i,j}^{\sigma_k} = \cases
A_{i-1,j}&\text{if $k=i-1$,}\\
A_{i+1,j}^{A_{i,i+1}}&\text{if $k=i<j-1$,}\\
A_{i,j-1}&\text{if $k=j-1>i$,}\\
A_{i,j+1}^{A_{j,j+1}}&\text{if $k=j$,}\\
A_{i,j}&\text{otherwise,}\\
\endcases
\qquad
A_{i,j}^{\sigma_k^{-1}} = \cases
A_{i-1,j}^{A_{i-1,i}^{-1}}&\text{if $k=i-1$,}\\
A_{i+1,j}&\text{if $k=i<j-1$,}\\
A_{i,j-1}^{A_{j-1,j}^{-1}}&\text{if $k=j-1>i$,}\\
A_{i,j+1}&\text{if $k=j$,}\\
A_{i,j}&\text{otherwise.}\\
\endcases
\tag{5}
$$

We shall work mainly with a particular type of pure braids.
These ``twist braids'' are defined as follows.
Given an increasingly ordered set  $I=\{i_1,\dots,i_r\}$, let
$$A_I = (A_{i_1,i_2})(A_{i_1,i_3}A_{i_2,i_3})(A_{i_1,i_4}
A_{i_2,i_4}A_{i_3,i_4})\cdots(A_{i_1,i_r}\cdots A_{i_{r-1},i_r}).
\tag{6}$$
We extend this definition to sets which are not increasingly ordered
(such as the vertex sets $V_k$ in~5.1) by first ordering, then
proceeding as above.
The conjugation action of an arbitrary braid
$\b\in B_n$ on the twist braid $A_I\in P_n$ takes the form
$$A_I^\b = A_{\omega(I)}^C,\tag{7}$$
where $\omega=\tau_n(\b)$, and $C=C(I,\b)$ is a pure braid that
may be computed recursively from~\thetag{5}.

\subhead 5.3 Braid monodromy
\endsubhead
We now extract the braid monodromy of $\A$ from an associated
braided wiring diagram $\W$.  By Theorem~3.4, the image of the
braid monodromy is contained in the pure braid group $P_n$.
We shall express the braid monodromy generators,
$\lambda_k:=\lambda(x_k)$, in terms of the standard
generators $A_{i,j}$.

The vertex set $V_k=\{i_1,\dots,i_r\}$ gives rise to a partition
$\P_k = L_k \cup V_k \cup U_k$, where $L_k$ (resp.~$U_k$) consists
of the indices of the wires below (resp.~above) the vertex $v_k$.
Let $I_k = \{j,j+1,\dots,j+r-1\}$ denote the {\it local index} of
$V_k$, where $j = |L_k| + 1$.  The local monodromy $\gamma_k$
around the point $y_k\in \Y$ is a full twist on $I_k$
given by the pure braid $A_{I_k}$ (compare~4.6).
Note that $A_{I_k}=\mu_{I_k}^2$, where
$$\mu_{I_k} =
(\s_{j}\cdots\s_{j+r-2})(\s_{j}\cdots\s_{j+r-3})\cdots
(\s_j\cdot\s_{j+1})\cdot(\s_j) \tag{8}$$
is a permutation braid---a half twist on $I_k$.
Also notice that the monodromy along a path from
$y_k'$ to $y_k''$ above (or below) the point $y_k$ is given by $\mu_{I_k}$.

To specify the braid monodromy of $\A$, it remains to identify the
conjugating braids $\b_k$ of~4.4, {\it Step 2}.  Choosing the paths
$\eta_k$ to coincide with $\xi$ between $y_j''$ and $y_{j+1}'$
for $j<k$, these conjugating braids may be expressed as
$\b_1 = \b_{0,1} = 1$, and
$\b_{k+1} = \b_{k,k+1}\cdot\mu_{I_k} \cdot \b_k$ for $k\ge 1$.
Hence the braid monodromy generators are given by
$$\lambda_k = A_{I_k}^{\b_k}. \tag{9}$$

Note that the state $\Pi_k$ of the braided wiring
diagram $\W$ is the image of the braid $\b_k$ in the symmetric group,
$\Pi_k=\tau_n(\b_k)$.  Note also that the vertex set $V_k$ and its local
index $I_k$ are related by $V_k=\Pi_k(I_k)$.  Thus, the braid
monodromy generators may be expressed solely in terms of pure braids:
$$\lambda_k =  A_{V_k}^{C_k}, \tag{10} $$
for certain $C_k\in P_n$.

\subhead 5.4  Conjugate arrangements
\endsubhead
Let $\A$ be an arrangement of $n$ lines in $\C^2$, with associated
braided wiring diagram $\W$ corresponding to a generic projection
$\pi:\C^2\to \C$ and admissible path $\xi$.  Let $\bar\A$ denote
the conjugate arrangement (see~3.8).  Clearly, the vertices of
$\bar\A$ are the complex conjugates of those of $\A$.  Thus,
$\pi$ is generic with respect to $\bar\A$, and $\bar\xi$ is admissible.
The corresponding braided diagram, $\overline\W$, is then obtained
from $\W$ by simply reversing the crossings of all the intermediate
braids.  Thus, the local indices of $\overline\W$ are given by
$\bar I_k=I_k$, the conjugating braids by
$\bar\b_{k+1} = \epsilon_n(\b_{k,k+1})\cdot\mu_{I_k} \cdot \bar\b_k$,
and the braid monodromy generators by
$\bar\lambda_k = A_{I_k}^{\bar\b_k}$.  From the proof of Theorem~3.9,
we see that the braid monodromy generators of the two conjugate
arrangements are related by
$$\bar\lambda_k=\epsilon_n\left( (\lambda_k^{-1})^
{\lambda_{k-1}^{-1}\cdots \lambda_1^{-1}} \right).$$

\subhead 5.5 Real arrangements
\endsubhead
If $\A$ is a real arrangement in $\C^2$ (that is, $\A$ is the
complexification of a line arrangement $\A_{\R}$ in $\R^2$),
then the defining polynomial $f(x,z)$ has real coefficients.
Consequently, the vertices of $\A$ all have real
coordinates, and their images under first-coordinate projection all
lie on the real axis in $\C$.  In this instance, we may take the path
$\xi:[0,1]\to\C$ to be a directed line segment along the real axis.
The resulting diagram $\W$ is unbraided---all the intermediate
braids $\b_{k,k+1}$ are trivial.  In other words, the diagram is a
wiring diagram in the combinatorial sense \cite{Go}, affine if $\A$
contains parallel lines (see also \cite{BLSWZ}).  In this instance,
the description of the braid monodromy given in~5.3 specializes to the
algorithm of Hironaka \cite{Hir}, modulo some notational differences.

Another description of the braid monodromy of an abstract (unbraided)
wiring diagram was provided by Cordovil and Fachada in \cite{CF} (see
also \cite{Cor}).  This description, based on Salvetti's work \cite{S1},
\cite{S2}, may be paraphrased as follows.

Recall that the vertex set $V_k=\{i_1,\dots,i_r\}$ gives rise to a
partition $\P_k = L_k \cup V_k \cup U_k$.  Let $\bar V_k = \{i\mid
i_1 \le i \le i_r\}$, and set $J_k = (\bar V_k \setminus V_k)\cap U_k$.
Let $B_{J_k} = \prod A_{j,i}$,
where the product is over all $j\in J_k$ and
$i\in V_k$ with $j<i$, taken in the natural order (so that $B_{J_k}$
is a subword of $\Delta^2=A_{1,\dots,n}$, equal to $1$ if $J_k=\emptyset$).
Then the braid monodromy generators may be expressed as
$\hat\l_k = A_{V_k}^{J_k}:= B_{J_k}^{-1}A_{V_k}B_{J_k}^{}$,
where $A_{V_k}$ is as defined in \thetag{6}.

Using the Artin representation, one can show that the braids $\l_k$
and $\hat\l_k$ are equal.  The action of the braid $\hat\l_k$ is given
in formulas \thetag{12} and \thetag{13} in section~6.  The same
formulas hold for $\l_k$, but the computation is more involved.
Thus, the two descriptions of the braid monodromy of a real arrangement
(or more generally, of an arbitrary wiring diagram) coincide.

\subhead 5.6 Markov moves
\endsubhead
For an arbitrary complex arrangement $\A$,
changes in the choices made in the construction of the braid
monodromy (see~3.5) give rise to changes
in the braided wiring diagram $\W$ associated to $\A$.
For instance, changing the basepoint $y_0$ may alter the initial braid
$\b_{0,1}$, while changes in the generic projection may alter
the order of the vertices.

We refer to these (and other) changes in a braided wiring diagram as
``Markov moves.''  In practice, these moves may be used to simplify a
braided wiring diagram associated to an arrangement $\A$ (and
consequently the braid monodromy generators of $\A$, as well).  We now
record these simplifying Markov moves and their effects on the braid
monodromy.  In the following, we record only the local index of a
vertex, so ``vertex $\{j,\dots,k\}$'' means ``a vertex with local index
$\{j,\dots,k\}$.''  Recall that, while we depict braided wiring diagrams
right to left, their intermediate braids are read left to right.

\remark{Geometric moves}
\item{(1)} Insert an arbitrary braid $\b_0$ at the beginning of the
braided wiring diagram.
\item{(2)} Insert an arbitrary braid $\b_{s+1}$ at the end of the
braided wiring diagram.
\item{(3)} Replace vertex $\{i,\dots,j\}$, then vertex
$\{k,\dots,l\}$ with
\itemitem{(a)} vertex $\{k,\dots,l\}$, then vertex
$\{i,\dots,j\}$, if $j<k$ or $i>l$.
\itemitem{(b)} braid $(\s_{k}\cdots \s_{i+1})
\cdot (\s_{k+1}\cdots \s_{i+2})\cdots
(\s_{l-1}\cdots \s_{j})$,
\itemitem{ }
then vertex $\{i,\dots,i+l-k\}$, then vertex $\{i+l-k,\dots,l\}$,
\itemitem{ }
then braid $(\s_{k-1}^{- 1}\cdots\s_i^{- 1})
\cdot(\s_{k}^{- 1}\cdots \s_{i+1}^{- 1})\cdots
(\s_{l-2}^{- 1}\cdots\s_{j-1}^{- 1})$, if $j=k$.
\itemitem{(c)} braid $(\s_{i-1}\cdots \s_{k})
\cdot (\s_{i}\cdots \s_{k+1})\cdots
(\s_{j-2}\cdots \s_{l-1})$,
\itemitem{ }
then vertex $\{j+k-l,\dots,j\}$, then vertex $\{k,\dots,j+k-l\}$,
\itemitem{ }
then braid $(\s_{i}^{- 1}\cdots\s_{k+1}^{- 1})
\cdot(\s_{i+1}^{- 1}\cdots \s_{k+2}^{- 1})\cdots
(\s_{j-1}^{- 1}\cdots\s_{l}^{- 1})$, if $i=l$.
\smallskip

\centerline{
\beginsmallgraph
\label {$k$} at (24.5, 1)
\label {$\vdots$} at (24.5, 2.2)
\label {$l$} at (24.5, 3)
\label {$\vdots$} at (24.5, 4.2)
\label {$j$} at (24.5, 5)
\edge from (24, 1) to (23, 1)
\edge from (24, 2) to (23.8, 2)
\edge from (23.8, 2) to (23.6, 2.33)
\edge from (23.4, 2.67) to (23.2, 3)
\edge from (23.2, 3) to (23, 3)
\edge from (24, 3) to (23.8, 3)
\edge from (23.8, 3) to (23.2, 2)
\edge from (23.2, 2) to (23, 2)
\edge from (24, 4) to (23, 4)
\edge from (24, 5) to (23, 5)
\edge from (23, 1) to (22.8, 1)
\edge from (22.8, 1) to (22.6, 1.33)
\edge from (22.4, 1.67) to (22.2, 2)
\edge from (22.2, 2) to (22, 2)
\edge from (23, 2) to (22.8, 2)
\edge from (22.8, 2) to (22.2, 1)
\edge from (22.2, 1) to (22, 1)
\edge from (23, 3) to (22.8, 3)
\edge from (22.8, 3) to (22.6, 3.33)
\edge from (22.4, 3.67) to (22.2, 4)
\edge from (22.2, 4) to (22, 4)
\edge from (23, 4) to (22.8, 4)
\edge from (22.8, 4) to (22.2, 3)
\edge from (22.2, 3) to (22, 3)
\edge from (23, 5) to (22, 5)
\edge from (22, 1) to (21, 1)
\edge from (22, 2) to (21.8, 2)
\edge from (21.8, 2) to (21.6, 2.33)
\edge from (21.4, 2.67) to (21.2, 3)
\edge from (21.2, 3) to (21, 3)
\edge from (22, 3) to (21.8, 3)
\edge from (21.8, 3) to (21.2, 2)
\edge from (21.2, 2) to (21, 2)
\edge from (22, 4) to (21, 4)
\edge from (22, 5) to (21, 5)
\edge from (21, 1) to (20, 1)
\edge from (21, 2) to (20, 2)
\edge from (21, 3) to (20.8, 3)
\edge from (20.8, 3) to (20.2, 5)
\edge from (20.2, 5) to (20, 5)
\edge from (21, 4) to (20.8, 4)
\edge from (20.8, 4) to (20.2, 4)
\edge from (20.2, 4) to (20, 4)
\edge from (21, 5) to (20.8, 5)
\edge from (20.8, 5) to (20.2, 3)
\edge from (20.2, 3) to (20, 3)
\edge from (20, 1) to (19.8, 1)
\edge from (19.8, 1) to (19.2, 3)
\edge from (19.2, 3) to (19, 3)
\edge from (20, 2) to (19.8, 2)
\edge from (19.8, 2) to (19.2, 2)
\edge from (19.2, 2) to (19, 2)
\edge from (20, 3) to (19.8, 3)
\edge from (19.8, 3) to (19.2, 1)
\edge from (19.2, 1) to (19, 1)
\edge from (20, 4) to (19, 4)
\edge from (20, 5) to (19, 5)
\edge from (19, 1) to (18, 1)
\edge from (19, 2) to (18, 2)
\edge from (19, 3) to (18.8, 3)
\edge from (18.8, 3) to (18.2, 4)
\edge from (18.2, 4) to (18, 4)
\edge from (19, 4) to (18.8, 4)
\edge from (18.8, 4) to (18.6, 3.67)
\edge from (18.4, 3.33) to (18.2, 3)
\edge from (18.2, 3) to (18, 3)
\edge from (19, 5) to (18, 5)
\edge from (18, 1) to (17, 1)
\edge from (18, 2) to (17.8, 2)
\edge from (17.8, 2) to (17.2, 3)
\edge from (17.2, 3) to (17, 3)
\edge from (18, 3) to (17.8, 3)
\edge from (17.8, 3) to (17.6, 2.67)
\edge from (17.4, 2.33) to (17.2, 2)
\edge from (17.2, 2) to (17, 2)
\edge from (18, 4) to (17.8, 4)
\edge from (17.8, 4) to (17.2, 5)
\edge from (17.2, 5) to (17, 5)
\edge from (18, 5) to (17.8, 5)
\edge from (17.8, 5) to (17.6, 4.67)
\edge from (17.4, 4.33) to (17.2, 4)
\edge from (17.2, 4) to (17, 4)
\edge from (17, 1) to (16, 1)
\edge from (17, 2) to (16, 2)
\edge from (17, 3) to (16.8, 3)
\edge from (16.8, 3) to (16.2, 4)
\edge from (16.2, 4) to (16, 4)
\edge from (17, 4) to (16.8, 4)
\edge from (16.8, 4) to (16.6, 3.67)
\edge from (16.4, 3.33) to (16.2, 3)
\edge from (16.2, 3) to (16, 3)
\edge from (17, 5) to (16, 5)
\label {$\sim$} at (15.5, 3)
\edge from (15, 1) to (14, 1)
\edge from (15, 2) to (14, 2)
\edge from (15, 3) to (14.8, 3)
\edge from (14.8, 3) to (14.2, 5)
\edge from (14.2, 5) to (14, 5)
\edge from (15, 4) to (14.8, 4)
\edge from (14.8, 4) to (14.2, 4)
\edge from (14.2, 4) to (14, 4)
\edge from (15, 5) to (14.8, 5)
\edge from (14.8, 5) to (14.2, 3)
\edge from (14.2, 3) to (14, 3)
\edge from (14, 1) to (13.8, 1)
\edge from (13.8, 1) to (13.2, 3)
\edge from (13.2, 3) to (13, 3)
\edge from (14, 2) to (13.8, 2)
\edge from (13.8, 2) to (13.2, 2)
\edge from (13.2, 2) to (13, 2)
\edge from (14, 3) to (13.8, 3)
\edge from (13.8, 3) to (13.2, 1)
\edge from (13.2, 1) to (13, 1)
\edge from (14, 4) to (13, 4)
\edge from (14, 5) to (13, 5)
\label {$i$} at (11.5, 1)
\label {$\vdots$} at (11.5, 2.2)
\label {$j$} at (11.5, 3)
\label {$\vdots$} at (11.5, 4.2)
\label {$l$} at (11.5, 5)
\edge from (11, 1) to (10, 1)
\edge from (11, 2) to (10, 2)
\edge from (11, 3) to (10.8, 3)
\edge from (10.8, 3) to (10.6, 3.33)
\edge from (10.4, 3.67) to (10.2, 4)
\edge from (10.2, 4) to (10, 4)
\edge from (11, 4) to (10.8, 4)
\edge from (10.8, 4) to (10.2, 3)
\edge from (10.2, 3) to (10, 3)
\edge from (11, 5) to (10, 5)
\edge from (10, 1) to (9, 1)
\edge from (10, 2) to (9.8, 2)
\edge from (9.8, 2) to (9.6, 2.33)
\edge from (9.4, 2.67) to (9.2, 3)
\edge from (9.2, 3) to (9, 3)
\edge from (10, 3) to (9.8, 3)
\edge from (9.8, 3) to (9.2, 2)
\edge from (9.2, 2) to (9, 2)
\edge from (10, 4) to (9.8, 4)
\edge from (9.8, 4) to (9.6, 4.33)
\edge from (9.4, 4.67) to (9.2, 5)
\edge from (9.2, 5) to (9, 5)
\edge from (10, 5) to (9.8, 5)
\edge from (9.8, 5) to (9.2, 4)
\edge from (9.2, 4) to (9, 4)
\edge from (9, 1) to (8, 1)
\edge from (9, 2) to (8, 2)
\edge from (9, 3) to (8.8, 3)
\edge from (8.8, 3) to (8.6, 3.33)
\edge from (8.4, 3.67) to (8.2, 4)
\edge from (8.2, 4) to (8, 4)
\edge from (9, 4) to (8.8, 4)
\edge from (8.8, 4) to (8.2, 3)
\edge from (8.2, 3) to (8, 3)
\edge from (9, 5) to (8, 5)
\edge from (8, 1) to (7.8, 1)
\edge from (7.8, 1) to (7.2, 3)
\edge from (7.2, 3) to (7, 3)
\edge from (8, 2) to (7.8, 2)
\edge from (7.8, 2) to (7.2, 2)
\edge from (7.2, 2) to (7, 2)
\edge from (8, 3) to (7.8, 3)
\edge from (7.8, 3) to (7.2, 1)
\edge from (7.2, 1) to (7, 1)
\edge from (8, 4) to (7, 4)
\edge from (8, 5) to (7, 5)
\edge from (7, 1) to (6, 1)
\edge from (7, 2) to (6, 2)
\edge from (7, 3) to (6.8, 3)
\edge from (6.8, 3) to (6.2, 5)
\edge from (6.2, 5) to (6, 5)
\edge from (7, 4) to (6.8, 4)
\edge from (6.8, 4) to (6.2, 4)
\edge from (6.2, 4) to (6, 4)
\edge from (7, 5) to (6.8, 5)
\edge from (6.8, 5) to (6.2, 3)
\edge from (6.2, 3) to (6, 3)
\edge from (6, 1) to (5, 1)
\edge from (6, 2) to (5.8, 2)
\edge from (5.8, 2) to (5.2, 3)
\edge from (5.2, 3) to (5, 3)
\edge from (6, 3) to (5.8, 3)
\edge from (5.8, 3) to (5.6, 2.67)
\edge from (5.4, 2.33) to (5.2, 2)
\edge from (5.2, 2) to (5, 2)
\edge from (6, 4) to (5, 4)
\edge from (6, 5) to (5, 5)
\edge from (5, 1) to (4.8, 1)
\edge from (4.8, 1) to (4.2, 2)
\edge from (4.2, 2) to (4, 2)
\edge from (5, 2) to (4.8, 2)
\edge from (4.8, 2) to (4.6, 1.67)
\edge from (4.4, 1.33) to (4.2, 1)
\edge from (4.2, 1) to (4, 1)
\edge from (5, 3) to (4.8, 3)
\edge from (4.8, 3) to (4.2, 4)
\edge from (4.2, 4) to (4, 4)
\edge from (5, 4) to (4.8, 4)
\edge from (4.8, 4) to (4.6, 3.67)
\edge from (4.4, 3.33) to (4.2, 3)
\edge from (4.2, 3) to (4, 3)
\edge from (5, 5) to (4, 5)
\edge from (4, 1) to (3, 1)
\edge from (4, 2) to (3.8, 2)
\edge from (3.8, 2) to (3.2, 3)
\edge from (3.2, 3) to (3, 3)
\edge from (4, 3) to (3.8, 3)
\edge from (3.8, 3) to (3.6, 2.67)
\edge from (3.4, 2.33) to (3.2, 2)
\edge from (3.2, 2) to (3, 2)
\edge from (4, 4) to (3, 4)
\edge from (4, 5) to (3, 5)
\label {$\sim$} at (2.5, 3)
\edge from (2, 1) to (1.8, 1)
\edge from (1.8, 1) to (1.2, 3)
\edge from (1.2, 3) to (1, 3)
\edge from (2, 2) to (1.8, 2)
\edge from (1.8, 2) to (1.2, 2)
\edge from (1.2, 2) to (1, 2)
\edge from (2, 3) to (1.8, 3)
\edge from (1.8, 3) to (1.2, 1)
\edge from (1.2, 1) to (1, 1)
\edge from (2, 4) to (1, 4)
\edge from (2, 5) to (1, 5)
\edge from (1, 1) to (0, 1)
\edge from (1, 2) to (0, 2)
\edge from (1, 3) to (0.8, 3)
\edge from (0.8, 3) to (0.2, 5)
\edge from (0.2, 5) to (0, 5)
\edge from (1, 4) to (0.8, 4)
\edge from (0.8, 4) to (0.2, 4)
\edge from (0.2, 4) to (0, 4)
\edge from (1, 5) to (0.8, 5)
\edge from (0.8, 5) to (0.2, 3)
\edge from (0.2, 3) to (0, 3)
\endgraph
}

\nopagebreak
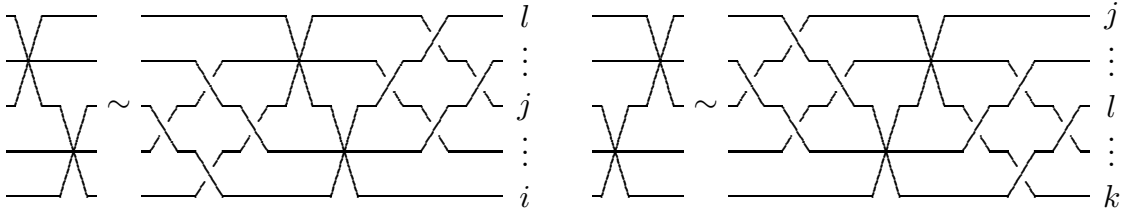
\centerline{Figure 2.  Moves 3b (left) and 3c (right)}
\endremark
\medskip
\remark{Further moves}
\item{(4)} Reduce the intermediate braid $\b_{k,k+1}$.
\item{(5)} Replace braid $\s_i$, then vertex
$\{j,\dots,k\}$ with
\itemitem{(a)} vertex $\{j,\dots,k\}$, then braid
$\s_i$, if $i<j-1$ or $i>k$.
\itemitem{(b)} braid $\s_j^{- 1}\cdots \s_{k-1}^{- 1}$, then
vertex $\{j+1,\dots,k+1\}$,  then braid
$\s_{k}\cdots \s_{j}$, if $i=k$.
\itemitem{(c)} braid $\s_{k-1}^{- 1}\cdots \s_{j}^{- 1}$, then
vertex $\{j-1,\dots,k-1\}$,  then braid
$\s_{j-1}\cdots \s_{k-1}$, if $i=j-1$.
\itemitem{(d)} vertex $\{j,\dots,k\}$, then braid
$\s_{j+k-i-1}$, if $j\le i \le k-1$.
\smallskip

\centerline{
\beginsmallgraph
\label {$j$} at (27.5, 1)
\label {$\vdots$} at (27.5, 2.2)
\label {$k$} at (27.5, 3)
\edge from (27, 1) to (26.8, 1)
\edge from (26.8, 1) to (26.2, 3)
\edge from (26.2, 3) to (26, 3)
\edge from (27, 2) to (26.8, 2)
\edge from (26.8, 2) to (26.2, 2)
\edge from (26.2, 2) to (26, 2)
\edge from (27, 3) to (26.8, 3)
\edge from (26.8, 3) to (26.2, 1)
\edge from (26.2, 1) to (26, 1)
\edge from (26, 1) to (25, 1)
\edge from (26, 2) to (25.8, 2)
\edge from (25.8, 2) to (25.6, 2.33)
\edge from (25.4, 2.67) to (25.2, 3)
\edge from (25.2, 3) to (25, 3)
\edge from (26, 3) to (25.8, 3)
\edge from (25.8, 3) to (25.2, 2)
\edge from (25.2, 2) to (25, 2)
\label {$\sim$} at (24.5, 2)
\edge from (24, 1) to (23.8, 1)
\edge from (23.8, 1) to (23.6, 1.33)
\edge from (23.4, 1.67) to (23.2, 2)
\edge from (23.2, 2) to (23, 2)
\edge from (24, 2) to (23.8, 2)
\edge from (23.8, 2) to (23.2, 1)
\edge from (23.2, 1) to (23, 1)
\edge from (24, 3) to (23, 3)
\edge from (23, 1) to (22.8, 1)
\edge from (22.8, 1) to (22.2, 3)
\edge from (22.2, 3) to (22, 3)
\edge from (23, 2) to (22.8, 2)
\edge from (22.8, 2) to (22.2, 2)
\edge from (22.2, 2) to (22, 2)
\edge from (23, 3) to (22.8, 3)
\edge from (22.8, 3) to (22.2, 1)
\edge from (22.2, 1) to (22, 1)
\label {$i$} at (20.5, 1)
\label {$j$} at (20.5, 2)
\label {$\vdots$} at (20.5, 3.2)
\label {$k$} at (20.5, 4)
\edge from (20, 1) to (19, 1)
\edge from (20, 2) to (19.8, 2)
\edge from (19.8, 2) to (19.2, 3)
\edge from (19.2, 3) to (19, 3)
\edge from (20, 3) to (19.8, 3)
\edge from (19.8, 3) to (19.6, 2.67)
\edge from (19.4, 2.33) to (19.2, 2)
\edge from (19.2, 2) to (19, 2)
\edge from (20, 4) to (19, 4)
\edge from (19, 1) to (18, 1)
\edge from (19, 2) to (18, 2)
\edge from (19, 3) to (18.8, 3)
\edge from (18.8, 3) to (18.2, 4)
\edge from (18.2, 4) to (18, 4)
\edge from (19, 4) to (18.8, 4)
\edge from (18.8, 4) to (18.6, 3.67)
\edge from (18.4, 3.33) to (18.2, 3)
\edge from (18.2, 3) to (18, 3)
\edge from (18, 1) to (17.8, 1)
\edge from (17.8, 1) to (17.2, 3)
\edge from (17.2, 3) to (17, 3)
\edge from (18, 2) to (17.8, 2)
\edge from (17.8, 2) to (17.2, 2)
\edge from (17.2, 2) to (17, 2)
\edge from (18, 3) to (17.8, 3)
\edge from (17.8, 3) to (17.2, 1)
\edge from (17.2, 1) to (17, 1)
\edge from (18, 4) to (17, 4)
\edge from (17, 1) to (16, 1)
\edge from (17, 2) to (16, 2)
\edge from (17, 3) to (16.8, 3)
\edge from (16.8, 3) to (16.6, 3.33)
\edge from (16.4, 3.67) to (16.2, 4)
\edge from (16.2, 4) to (16, 4)
\edge from (17, 4) to (16.8, 4)
\edge from (16.8, 4) to (16.2, 3)
\edge from (16.2, 3) to (16, 3)
\edge from (16, 1) to (15, 1)
\edge from (16, 2) to (15.8, 2)
\edge from (15.8, 2) to (15.6, 2.33)
\edge from (15.4, 2.67) to (15.2, 3)
\edge from (15.2, 3) to (15, 3)
\edge from (16, 3) to (15.8, 3)
\edge from (15.8, 3) to (15.2, 2)
\edge from (15.2, 2) to (15, 2)
\edge from (16, 4) to (15, 4)
\edge from (15, 1) to (14.8, 1)
\edge from (14.8, 1) to (14.6, 1.33)
\edge from (14.4, 1.67) to (14.2, 2)
\edge from (14.2, 2) to (14, 2)
\edge from (15, 2) to (14.8, 2)
\edge from (14.8, 2) to (14.2, 1)
\edge from (14.2, 1) to (14, 1)
\edge from (15, 3) to (14, 3)
\edge from (15, 4) to (14, 4)
\label {$\sim$} at (13.5, 2.5)
\edge from (13, 1) to (12.8, 1)
\edge from (12.8, 1) to (12.6, 1.33)
\edge from (12.4, 1.67) to (12.2, 2)
\edge from (12.2, 2) to (12, 2)
\edge from (13, 2) to (12.8, 2)
\edge from (12.8, 2) to (12.2, 1)
\edge from (12.2, 1) to (12, 1)
\edge from (13, 3) to (12, 3)
\edge from (13, 4) to (12, 4)
\edge from (12, 1) to (11, 1)
\edge from (12, 2) to (11.8, 2)
\edge from (11.8, 2) to (11.2, 4)
\edge from (11.2, 4) to (11, 4)
\edge from (12, 3) to (11.8, 3)
\edge from (11.8, 3) to (11.2, 3)
\edge from (11.2, 3) to (11, 3)
\edge from (12, 4) to (11.8, 4)
\edge from (11.8, 4) to (11.2, 2)
\edge from (11.2, 2) to (11, 2)
\label {$j$} at (9.5, 1)
\label {$\vdots$} at (9.5, 2.2)
\label {$k$} at (9.5, 3)
\edge from (9, 1) to (8, 1)
\edge from (9, 2) to (8.8, 2)
\edge from (8.8, 2) to (8.2, 3)
\edge from (8.2, 3) to (8, 3)
\edge from (9, 3) to (8.8, 3)
\edge from (8.8, 3) to (8.6, 2.67)
\edge from (8.4, 2.33) to (8.2, 2)
\edge from (8.2, 2) to (8, 2)
\edge from (9, 4) to (8, 4)
\edge from (8, 1) to (7.8, 1)
\edge from (7.8, 1) to (7.2, 2)
\edge from (7.2, 2) to (7, 2)
\edge from (8, 2) to (7.8, 2)
\edge from (7.8, 2) to (7.6, 1.67)
\edge from (7.4, 1.33) to (7.2, 1)
\edge from (7.2, 1) to (7, 1)
\edge from (8, 3) to (7, 3)
\edge from (8, 4) to (7, 4)
\edge from (7, 1) to (6, 1)
\edge from (7, 2) to (6.8, 2)
\edge from (6.8, 2) to (6.2, 4)
\edge from (6.2, 4) to (6, 4)
\edge from (7, 3) to (6.8, 3)
\edge from (6.8, 3) to (6.2, 3)
\edge from (6.2, 3) to (6, 3)
\edge from (7, 4) to (6.8, 4)
\edge from (6.8, 4) to (6.2, 2)
\edge from (6.2, 2) to (6, 2)
\edge from (6, 1) to (5.8, 1)
\edge from (5.8, 1) to (5.6, 1.33)
\edge from (5.4, 1.67) to (5.2, 2)
\edge from (5.2, 2) to (5, 2)
\edge from (6, 2) to (5.8, 2)
\edge from (5.8, 2) to (5.2, 1)
\edge from (5.2, 1) to (5, 1)
\edge from (6, 3) to (5, 3)
\edge from (6, 4) to (5, 4)
\edge from (5, 1) to (4, 1)
\edge from (5, 2) to (4.8, 2)
\edge from (4.8, 2) to (4.6, 2.33)
\edge from (4.4, 2.67) to (4.2, 3)
\edge from (4.2, 3) to (4, 3)
\edge from (5, 3) to (4.8, 3)
\edge from (4.8, 3) to (4.2, 2)
\edge from (4.2, 2) to (4, 2)
\edge from (5, 4) to (4, 4)
\edge from (4, 1) to (3, 1)
\edge from (4, 2) to (3, 2)
\edge from (4, 3) to (3.8, 3)
\edge from (3.8, 3) to (3.6, 3.33)
\edge from (3.4, 3.67) to (3.2, 4)
\edge from (3.2, 4) to (3, 4)
\edge from (4, 4) to (3.8, 4)
\edge from (3.8, 4) to (3.2, 3)
\edge from (3.2, 3) to (3, 3)
\label {$\sim$} at (2.5, 2.5)
\edge from (2, 1) to (1, 1)
\edge from (2, 2) to (1, 2)
\edge from (2, 3) to (1.8, 3)
\edge from (1.8, 3) to (1.6, 3.33)
\edge from (1.4, 3.67) to (1.2, 4)
\edge from (1.2, 4) to (1, 4)
\edge from (2, 4) to (1.8, 4)
\edge from (1.8, 4) to (1.2, 3)
\edge from (1.2, 3) to (1, 3)
\edge from (1, 1) to (0.8, 1)
\edge from (0.8, 1) to (0.2, 3)
\edge from (0.2, 3) to (0, 3)
\edge from (1, 2) to (0.8, 2)
\edge from (0.8, 2) to (0.2, 2)
\edge from (0.2, 2) to (0, 2)
\edge from (1, 3) to (0.8, 3)
\edge from (0.8, 3) to (0.2, 1)
\edge from (0.2, 1) to (0, 1)
\edge from (1, 4) to (0, 4)
\endgraph
}

\nopagebreak
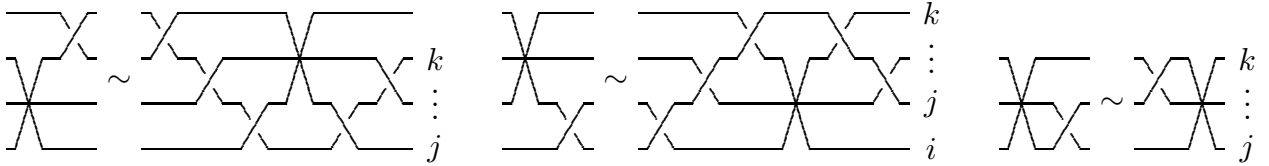
\centerline{Figure 3.  Moves 5b, 5c, and 5d (left to right)}
\endremark
\medskip
The parity of the braids in move (3) and move (5) may be switched.
For instance, one can replace braid $\s_i^{-1}$, then vertex
$\{j,\dots,k\}$, with braid $\s_j\cdots \s_{k-1}$, then
vertex $\{j+1,\dots,k+1\}$,  then braid
$\s_{k}^{- 1}\cdots \s_{j}^{- 1}$
if $i=k$ (move 5b).

Note that the triangle-switches and flips discussed in \cite{BLSWZ}
and \cite{Fa}
in the context of (unbraided) wiring diagrams may be accomplished
using Markov moves of types (3), (4), and (5).

\proclaim{Theorem 5.7} The braid monodromy of a braided wiring diagram
is invariant under Markov moves:  If the braided wiring diagram
$\widehat\W$ is obtained from the braided wiring diagram $\W$ by a
finite sequence of Markov moves of types (1)--(5) and their inverses,
then the braid monodromy homomorphisms $\lambda$ of $\W$ and $\hat\lambda$
of $\widehat\W$ are braid-equivalent.
\endproclaim
\demo{Sketch of Proof}
One can check, either algebraically or by drawing the appropriate
braids, that the (only) effects on the braid
monodromy of the Markov moves listed above are as follows:

\item{(1)} Global conjugation by $\b_0$:\quad
$\hat\l =\b_0^{-1}\cdot \l \cdot \b_0$.
\item{(2)} None.
\item{(3)} Suppose the vertices in question are the $k^{\text{th}}$
and $(k+1)^{\text{st}}$ vertices of $\W$ (resp.~$\widehat\W$).
Then the corresponding braid monodromy generators, $\l_k,
\l_{k+1}$ and $\hat\l_k, \hat\l_{k+1}$,
satisfy:
\itemitem{(a)} $\hat\l_k = \l_{k+1}$ and $\hat\l_{k+1}=\l_k$.
Note that the permutation braids $\mu_{I_k}$ and $\mu_{I_{k+1}}$ commute
in this instance. Thus we can write
$\hat\l_k = \l_k\cdot\l_{k+1}\cdot\l_k^{- 1}=\l_{k+1}$,
and $\hat\l = \l\circ \s_k$.
\itemitem{(b)} $\hat\l_k = \l_k\cdot\l_{k+1}\cdot\l_k^{- 1}$
and $\hat\l_{k+1}=\l_k$.  Thus $\hat\l = \l\circ \s_k$.
\itemitem{(c)}  $\hat\l_k = \l_{k+1}$ and
$\hat\l_{k+1}=\l_{k+1}^{- 1}\cdot\l_k\cdot\l_{k+1}$.
Thus $\hat\l=\l\circ\s_k^{-1}$.
\item{(4)} None.
\item{(5)} None.

\noindent Hence $\l$ and $\hat\l$ are braid-equivalent.
\quad\qed
\enddemo

Recall that curves in the same connected
component of an equisingular family have braid-equivalent
monodromies (see \cite{L3} and ~3.5).  In particular,
this is the case for arrangements which
are lattice-isotopic \cite{R2}.

\proclaim{Corollary 5.8} Let $\A$ and $\A'$ be lattice-isotopic
arrangements in $\C^2$, with associated braided wiring diagrams $\W$
and $\W'$. Then $\W'$ may be obtained from $\W$ via a finite sequence
of Markov moves and their inverses.
\endproclaim

\head 6. The group of a complex arrangement
\endhead

We now turn our attention the fundamental group $G$ of the complement
of the line arrangement $\A$ in $\C^2$.  In this section, we describe
the braid monodromy and Arvola presentations of $G$, show that
they are Tietze-I equivalent, and derive some homotopy type
consequences.

\subhead 6.1 Presentations
\endsubhead
Using the (pure) braid monodromy generators $\{\l_k\}$
from~\thetag{9} and the procedure described in~4.1, we obtain the braid
monodromy presentation
$$
G = \langle t_1,\dots,t_n\mid \l_k(t_i)=t_i\ \text{for}\ i
\in \check V_k\ \text{and}\ 1\le k\le s\rangle,$$
where $\check V_k = V_k \setminus \max V_k$.

We may also use the braided wiring diagram $\W$ to find the Arvola
presentation of $G$.  This presentation is obtained by applying the
Arvola algorithm \cite{Ar}, \cite{OT}, to $\W$.  Explicitly, we sweep a
vertical line across the braided wiring diagram from right to left,
introducing relations and keeping track of conjugations as we pass
through the vertices $v_k$ and the braids $\b_{k,k+1}$.  Note that the
braids $\b_{k,k+1}$ do not give rise to any relations, but do cause
additional conjugations.  It is convenient to express the generators
and these relations and conjugations in terms of the inverses of the
generators typically used in the Arvola algorithm.  Apart from this
notational difference, our description of this method follows Falk's
discussion in \cite{Fa} of the Randell algorithm for real arrangements
(to which the Arvola algorithm specializes).

Recall that the symbol $[g_1,\dots,g_r]$ denotes the family of $r-1$
relations
$$g_1\cdot g_2\cdots g_m = g_2\cdots g_m\cdot g_1 = \cdots\cdots =
g_m\cdot g_1\cdots g_{m-1}.$$
The Arvola presentation of the group $G$ is given by
$$
G = \langle t_1,\dots,t_n\mid \RR_1,\dots,\RR_s\rangle,$$
where if $V_k=\{i_1,\dots,i_r\}$, then $\RR_k$ denotes the family of
relations $[x_{i_1}(k),\dots,x_{i_r}(k)]$.  The word $x_i(k)$ denotes
the meridian about wire $i$ at state $\P_k$.  Let $y_i(k)$ denote the
meridian about wire $i$ between vertex $k$ and braid $k$.  Then we have
$$y_i(k)=\cases
x_i(k)&\text{if $i\notin V_k$,}\\
x_{i_1}(k)\cdots x_{i_l}(k)\cdot
x_i(k)\cdot(x_{i_1}(k)\cdots x_{i_l}(k))^{-1}&\text{if $i=i_{l+1}\in
V_k$.}
\endcases
\tag{11}$$
(see below).  The words $x_i(k)$ satisfy the recursion $x_i(1)=t_i$,
and $x_i(k+1)=\hat\b_{k,k+1}(y_i(k))$ for $k>0$, where $\hat\b_{k,k+1}$ records
the effect of the braiding $\b_{k,k+1}$ on the meridian $y_i(k)$ as indicated
below.

\bigskip

\centerline{
\beginsmallgraph
\edge from (0, 0) to (2, 3)
\edge from (0, 2) to (2, 1)
\edge from (0, 3) to (2, 0)
\label {$x_1$} at (2.5, 0)
\label {$x_2$} at (2.5, 1)
\label {$\vdots$} at (2.5, 2.3)
\label {$x_r$} at (2.5, 3)
\label {$y_1$} at (-.5, 3)
\label {$y_2$} at (-.5, 2)
\label {$\vdots$} at (-.5, 1.3)
\label {$y_r$} at (-.5, 0)
\edge from (10, 0) to (12, 2)
\edge from (10, 2) to (10.7, 1.3)
\edge from (12, 0) to (11.3, .7)
\label {$y_2$} at (12.5, 2)
\label {$y_1$} at (12.5, 0)
\label {$y_1\cdot y_2\cdot y_1^{-1}$} at (8.25, 0)
\label {$y_1$} at (9.5, 2)
\edge from (20, 2) to (22, 0)
\edge from (20, 0) to (20.7, .7)
\edge from (22, 2) to (21.3, 1.3)
\label {$y_2$} at (22.5, 2)
\label {$y_1$} at (22.5, 0)
\label {$y_2^{-1}\cdot y_1\cdot y_2$} at (18, 2)
\label {$y_2$} at (19.5, 0)
\label {$x_i$} at (2.5, -1.25)
\label {$\leftarrow$} at (1, -1.25)
\label {$\mu(x_i)$} at (-1, -1.25)
\label {$y_i$} at (12.5, -1.25)
\label {$\leftarrow$} at (11, -1.25)
\label {$\sigma_1(y_i)$} at (9, -1.25)
\label {$y_i$} at (22.5, -1.25)
\label {$\leftarrow$} at (21, -1.25)
\label {$\sigma_1^{-1}(y_i)$} at (19, -1.25)
\endgraph
}

\medskip
\nopagebreak
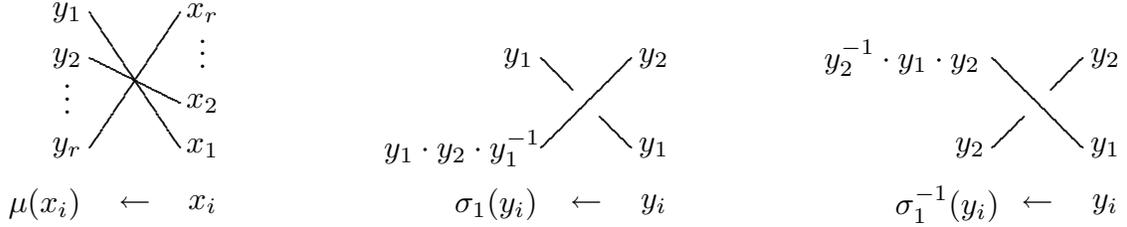
\centerline{Figure 4.  Conjugations in Arvola's algorithm}
\medskip
\noindent
Note that, locally, the conjugation arising from a vertex (resp.~braiding)
coincides with the action of a permutation braid (resp.~an
elementary braid or its inverse).  To make the description of
$\hat\b_{k,k+1}$ explicit, we require some notation.

The state $\P_k$ of the braided wiring diagram $\W$ is a permutation
of $\{1,\dots,n\}$, recording the relative heights of the wires at this state.
Recall the local index $I_k$ of the vertex set $V_k$ of $\W$,  and the
associated permutation braid $\mu_{I_k}$.  Let $\bar\mu_{I_k}=\tau_n(\mu_{I_k})$
denote the permutation induced by $\mu_{I_k}$, and
let $\bar\b_{k,k+1}=\tau_n(\b_{k,k+1})$.  It is easily
seen that $\P_{k+1}=\bar\b_{k,k+1}\cdot\bar\mu_{I_k}\cdot\P_k$.

Note that the sets $\{x_i(k)\}$ and $\{y_i(k)\}$ generate the free group
$F_n=\langle t_1,\dots,t_n\rangle$.  Define $\phi_k, \psi_k \in
\Aut(F_n)$ by
$$\phi_k(t_q)=x_i(k)\ \text{if}\ \P_k(q)=i\quad
\text{and}\quad
\psi_k(t_q)=x_i(k)\ \text{if}\ \bar\mu_{I_k}\cdot\P_k(q)=i.$$
If $\P_{k+1}(q)=i$, then the effect of the braiding $\b_{k,k+1}$ on
$y_i(k)$ may be expressed as $\hat\b_{k,k+1}(y_i(k)) =
\b_{k,k+1}\cdot\psi_k(t_q) = \psi_k(\b_{k,k+1}(t_q))$.

\proclaim{Lemma 6.2} We have $\psi_k=\mu_{I_k}\cdot\phi_k$.
\endproclaim
\demo{Proof} Compute:
$$\mu_{I_k}\cdot\phi_k(t_q)= 
\cases
\phi_k(t_q)&\text{if $q\notin I_k$,}\\
\phi_k(t_j\cdots t_{j+l-1}\cdot t_{j+l}\cdot (t_j\cdots t_{j+l-1})^{-1})&
\text{if $q=j+r-1-l \in I_k$.}
\endcases$$
Checking that this agrees with the description of the meridians
$y_i(k)$ given in \thetag{11}, we have
$\mu_{I_k}\cdot\phi_k(t_q)=y_i(k)=\psi_k(t_q)$ for all $q$.\quad\qed
\enddemo

We now show that the meridians $x_i(k)$ may be expressed in terms
of the conjugating braids $\b_k$ from the braid monodromy
constructions of~4.4 and 5.3.  Recall that these braids are defined by
$\b_1=1$, and $\b_{k+1}=\b_{k,k+1}\cdot\mu_{I_k}\cdot \beta_k$ for $k\ge 1$.

\proclaim{Proposition 6.3} If wire $i$ is at height $q$ at state
$\P_{k+1}$ in the braided wiring diagram $\W$ (that is,
$\P_{k+1}(q)=i$), then $x_i(k+1)=\b_{k+1}(t_q)$.
\endproclaim
\demo{Proof} We use induction on $k$, with the case $k=0$ trivial.

In general, assume $\P_{k+1}(q)=i$.  We have
$$x_i(k+1)=\hat\b_{k,k+1}(y_i(k))= 
\b_{k,k+1}\cdot\mu_{I_k}\cdot\phi_k(t_q) =
\b_{k,k+1}\cdot\mu_{I_k}\cdot \b_k(t_q) =
\b_{k+1}(t_q),$$
using the above lemma, the inductive hypothesis to identify
$\phi_k=\b_k$, and the identification of the braids $\b_k$ from~5.3.\quad\qed
\enddemo

We may now state and prove the main theorem of this section.

\proclaim{Theorem 6.4} The braid monodromy and Arvola presentations
of the group $G$ of the arrangement $\A$ are Tietze-I equivalent.
\endproclaim
\demo{Proof} Let $V=V_k=\{i_1,\dots,i_r\}$ denote the
$k^{\text{th}}$ vertex of a braided wiring diagram $\W$ associated
to $\A$, with local index $I=I_k=\{j,\dots,j+r-1\}$.  Write $\b=\b_k$ and
$\a=\a_k=A_I^\b$, and recall that $\check V = V \setminus \max V$.  Using
Proposition 6.3, we may express the family
$\RR_k$ of Arvola relations as
$$[\b(t_j),\dots,\b(t_{j+r-1})].$$
We will show that these $r-1$ relations and the braid monodromy relations
$\a(t_i)=t_i$, $i\in \check V$ are equivalent.  It is
easy to see
that the $r-1$ Arvola relations above are equivalent to
$\b(A_I(t_i))=\b(t_i)$, $i\in \check I$.

Using Proposition 6.3 again, we have $\b(t_i)=x_{i_p}(k)$ if
$i=j+p-1\in I$. Consequently, $\b(t_i)$ is some conjugate of $t_{i_p}$,
say $\b(t_i)=w_p\cdot t_{i_p}\cdot w_p^{-1}$.  A computation shows
that $\a(t_{i_p})=\a(w_p^{-1})\cdot \b(A_I(t_i))\cdot \a(w_p)$.
Therefore the braid monodromy relation $\a(t_{i_p})=t_{i_p}$ may be
expressed as $\b(A_I(t_i))=\a(w_p) \cdot t_{i_p} \cdot \a(w_p^{-1})$.
Now it follows from the relations $\a(t_i)=t_i$ for $i\in \check V$ that
$\a(t_i)=t_i$ for all $i$.  Thus $\a(w_p)=w_p$, and
the braid monodromy relation above is clearly equivalent to the
corresponding Arvola relation.\quad\qed
\enddemo

Since the braid monodromy and Arvola presentations of the group $G$
of $\A$ are Tietze-I equivalent, the associated 2-complexes are
homotopy equivalent.  Furthermore, Libgober's theorem \cite{L1}
stated in~4.2 provides a homotopy equivalence between the
braid monodromy presentation 2-complex and the complement
$\C^2\setminus\A$.  Thus we obtain the following corollary.

\proclaim{Corollary 6.5} The complement of a complex arrangement $\A$
in $\C^2$ has the homotopy type of the 2-complex modeled on the
Arvola presentation of the group $G$.
\endproclaim
Prior to our obtaining this result, Arvola informed us that he had
a proof of it.  We are not cognizant of the details of that proof.

\subhead 6.6 Real arrangements
\endsubhead
If $\A$ is a real arrangement in $\C^2$, then, as noted in 5.5, the
braided wiring diagram $\W$ is unbraided, so is a (possibly affine)
wiring diagram.
In this instance, Arvola's algorithm specializes to that of Randell
\cite{R1}.  Using Theorem 6.4, we obtain:

\proclaim{Corollary 6.7} The braid monodromy and Randell presentations
of the group $G$ of a real arrangement $\A$ in $\C^2$ are Tietze-I
equivalent.
\endproclaim

As above, combining this result with Libgober's theorem yields the
following corollary, which constitutes the main result of Falk~\cite{Fa}.

\proclaim{Corollary 6.8} The complement of a real arrangement
$\A$ in $\C^2$ has the homotopy type of the 2-complex modeled on the
Randell presentation of the group $G$.
\endproclaim

The braid monodromy and Randell presentations of $G$ may be obtained
immediately from the description of the generators $\hat\l_k=A_{V_k}^{J_k}$
in terms of pure braids provided by \cite{CF} and described in~5.5.
This is accomplished by finding the action of the braids $\hat\l_k$.  For
the sake of completeness, we find the action of $\hat\l_k$ on the
entire free group $F_n$.

Write $V=V_k$ and $J=J_k$.  If $V=\{i_1,\dots,i_r\}$, let
$t_V = t_{i_1}\cdots t_{i_r}$ (set $t_V=1$ if $V=\emptyset$).  For
$i\in \bar V\setminus V$, let $V^{<i}=\{i_1,\dots,i_q\}$ and
$V^{>i}=\{i_{q+1},\dots,i_r\}$ if $i_q<i<i_{q+1}$.  If
$J =\emptyset$, a straightforward computation yields
$$A_V(t_i)=\cases  t_V^{}\cdot t_i^{}\cdot t_V^{-1}
&\text{ if } i\in V,\\
[t_{V^{<i}},t_{V^{>i}}]\cdot t_i \cdot
[t_{V^{<i}},t_{V^{>i}}]^{-1}
&\text{ if } i\in \bar V \setminus V,\\
t_i &\text{ otherwise.}
\endcases \tag{12}$$
If $J\neq\emptyset$, let
$$z_{J,i}^{} = \cases 1&\text{if $i<j_1$ or $i \in J$ or $i>i_r$,}\\
t_{J^{<i}}&\text{if $i\in\bar J\setminus J$,}\\  t_{J}&\text{if
$\max J < i \le i_r$,}\\
\endcases$$
and define $\gamma_J \in \Aut(F_n)$ by $\gamma_J(t_i) =
z_{J,i}^{}\cdot t_i^{}\cdot z_{J,i}^{-1}$.
Induction on $|J|$, starting from \thetag{12}, yields:
$$A_V^{J}(t_i)=\cases
z_{J,i}^{-1}\cdot \gamma_J (t_V^{}\cdot t_i^{}
\cdot t_V^{-1}) \cdot z_{J,i}
&\text{if $i\in V$,}\\
z_{J,i}^{-1}\cdot \gamma_J ([t_{V^{<i}},t_{V^{>i}}]\cdot t_i
\cdot[t_{V^{<i}},t_{V^{>i}}]^{-1})\cdot z_{J,i}
&\text{if $i\in \bar V\setminus (V\cup J)$,}\\
t_i &\text{if $i \in J$ or $i \notin \bar V$.}
\endcases\tag{13}$$

\proclaim{Proposition 6.9} Let $\A$ be a real arrangement.
If $\W$ is an associated wiring diagram
with vertex sets $V_k$ and conjugating sets $J_k$, then
the braid monodromy and Randell presentations
of the group $G(\A)$ are given by
$$\align
G&=\langle t_1,\dots,t_n\mid \gamma_k(t_{V_k}\cdot t_i \cdot
t_{V_k}^{-1}) =
\gamma_k(t_i)\ \text{for}\ i
\in \check V_k\ \text{and}\ 1\le k\le s\rangle\\
&=\langle t_1,\dots,t_n\mid \RR_1,\dots,\RR_s\rangle,
\endalign$$
where $\gamma_k=\gamma_{J_k}$, and
$\RR_k$ denotes the family of relations
$[\gamma_k(t_{i_1}),\dots,\gamma_k(t_{i_r})]$.
\endproclaim

\head 7. Applications and examples
\endhead

In this section, we demonstrate the techniques described above by
means of several explicit examples and provide some applications.

\subhead 7.1 The intersection lattice
\endsubhead
Let $\A=\{ H_1,\dots,H_n\}$ be an arrangement,
and let $L(\A)$ be the ranked poset of non-empty intersections
of $\A$, ordered by reverse inclusion, and
with rank function given by codimension.  Two arrangements
$\A$ and $\A'$ are {\it lattice-isomorphic} if there
is an order-preserving bijection $\pi: L(\A)\to L(\A')$ (see
\cite{OT} for further details).

Let $\A$ be an arrangement of $n$ lines in $\C^2$ with $s$ vertices.
Choose (arbitrary) orderings of the lines and vertices of $\A$.
Then the intersection lattice may be encoded simply by a map
$V: \{1,\dots,s\} \to \S(n)$, where $\S(n)$ denotes
the set of all subsets of $\{1,\dots,n\}$, and
$V(k)=V_k$  is the $k^{\text{th}}$ vertex set.
Two arrangements $\A$ and $\A'$ in $\C^2$ are lattice-isomorphic
if, upon ordering their respective lines and vertices,
there exist permutations $\pi\in \Sigma_s$ and
$\rho\in \Sigma_n$ such that $V'_{\pi(k)}=\rho(V_k)$.

\proclaim{Theorem 7.2} Line arrangements with braid-equivalent monodromies
are lattice-isomor\-phic.
\endproclaim

\demo{Proof} First recall that, if $A_I$ is
one of the (extended) generators of $P_n$ specified in \thetag{6},
and $\beta\in B_n$ with $\tau_n(\beta)=\omega$, then
$A_I^{\beta}=A_{\omega(I)}^C$,  for some $C\in P_n$, see \thetag{7}.
Also note that, since the abelianization of $P_n$ is a free abelian group
on the images of the standard generators $A_{i,j}$, if a pure braid $\gamma$
can be written as $\gamma=B^{-1}\cdot A_I^{\pm 1}\cdot B$, for some
$B\in P_n$, then the indexing set $I$ is uniquely determined by
$\gamma$.

Now let $\A$ and $\A'$ be line arrangements with braid-equivalent
monodromies $\a$ and $\a'$.
Note that the braid monodromy construction determines orderings
of the lines and vertices of the arrangements.
Write $\a(x_k)=A_{V_k}^{C_k}$ and $\a'(x_k)=A_{V'_k}^{C'_k}$,
with $C_k, C'_k\in P_n$ as in \thetag{10}.
By assumption, $\a'\circ \psi = \conj_{\phi}\circ\a$,
with $\psi\in B_s$, $\phi\in B_n$.
Write $\psi(x_k)= z_k^{-1} \cdot x_{\pi(k)} \cdot z_k$,
where $\pi=\tau_s(\psi)$.  Also, set
$\rho=\tau_n(\phi)$.  The braid-equivalence then reads:
$$A_{V'_{\pi(k)}}^{C'_k\a'(z_k)} =
A_{\rho(V_k)}^{B_kC_k^{\phi}}.$$
Since both exponents in the above equation are pure braids,
we conclude that $V'_{\pi(k)}=\rho(V_k)$, as
required.  \qquad\qed
\enddemo

\subhead Example 7.3
\endsubhead
One can easily find pairs of arrangements whose groups are isomorphic,
yet whose monodromies are not equivalent.  For instance, consider
arrangements with defining polynomials $Q(\A)=xy(x-y)$ and $Q(\A')=xy(x-1)$,
respectively.  It is readily checked that $G(\A)=P_3$ is isomorphic
to $G(\A')=F_2\times F_1$.  Furthermore, it can be seen that
$\C^2\setminus \A=(S^3\setminus 3 \text{ Hopf circles})\times \R^+$
is diffeomorphic to $\C^2\setminus \A'=
S^1\times (S^2 \setminus 3 \text{ points})\times \R^+$.
On the other hand, the respective pure braid monodromies,
$\l:F_1\to P_3$, $\l(x)=A_{1,2,3}$, and $\l':F_2\to P_3$,
$\l'(x_1)=A_{1,2}$, $\l'(x_2)=A_{1,3}$, are obviously
not equivalent.  This may be explained by the fact that there
is no ambient diffeomorphism of $\C^2$ taking $\A$ to $\A'$.

While these examples do show that the braid monodromy of a plane curve
carries more information than the fundamental group of the complement,
they are unsatisfying for several reasons.
Combinatorially, $L(\A)$ is a lattice, while $L(\A')$ is merely a poset.
Geometrically, $\A$ is transverse to the line at infinity, while $\A'$ is not.
This being the case, these examples do not address Libgober's question
\cite{L3}, which was posed for plane curves that are transverse to the
line at infinity in $\CP^2$.  We now present examples which do fit
into this framework.

\subhead 7.4 Falk arrangements
\endsubhead
Consider the arrangements
of Falk \cite{Fa} defined by
$$Q=yz(x+y)(x-y)(x+z)(x-z)\quad\text{and}\quad
Q'=yz(x+z)(x-z)(y-z)(x-y-z).$$
Taking generic sections, we get a pair of real line arrangements
$\A$ and $\A'$ in $\C^2$ which are transverse to the line at
infinity, and have the same numbers of double and triple points.
Wiring diagrams for these line arrangements are depicted below.

\medskip

\centerline{
\beginsmallgraph
\label {\ } at (26.5, 0)
\label {$1$} at (25.5, 1)
\label {$2$} at (25.5, 2)
\label {$3$} at (25.5, 3)
\label {$4$} at (25.5, 4)
\label {$5$} at (25.5, 5)
\label {$6$} at (25.5, 6)
\edge from (25, 1) to (24.8, 1)
\edge from (24.8, 1) to (24.2, 3)
\edge from (24.2, 3) to (24, 3)
\edge from (25, 2) to (24.8, 2)
\edge from (24.8, 2) to (24.2, 2)
\edge from (24.2, 2) to (24, 2)
\edge from (25, 3) to (24.8, 3)
\edge from (24.8, 3) to (24.2, 1)
\edge from (24.2, 1) to (24, 1)
\edge from (25, 4) to (24, 4)
\edge from (25, 5) to (24, 5)
\edge from (25, 6) to (24, 6)
\edge from (24, 1) to (23, 1)
\edge from (24, 2) to (23, 2)
\edge from (24, 3) to (23.8, 3)
\edge from (23.8, 3) to (23.2, 5)
\edge from (23.2, 5) to (23, 5)
\edge from (24, 4) to (23.8, 4)
\edge from (23.8, 4) to (23.2, 4)
\edge from (23.2, 4) to (23, 4)
\edge from (24, 5) to (23.8, 5)
\edge from (23.8, 5) to (23.2, 3)
\edge from (23.2, 3) to (23, 3)
\edge from (24, 6) to (23, 6)
\edge from (23, 1) to (22, 1)
\edge from (23, 2) to (22.8, 2)
\edge from (22.8, 2) to (22.2, 3)
\edge from (22.2, 3) to (22, 3)
\edge from (23, 3) to (22.8, 3)
\edge from (22.8, 3) to (22.2, 2)
\edge from (22.2, 2) to (22, 2)
\edge from (23, 4) to (22, 4)
\edge from (23, 5) to (22, 5)
\edge from (23, 6) to (22, 6)
\edge from (22, 1) to (21.8, 1)
\edge from (21.8, 1) to (21.2, 2)
\edge from (21.2, 2) to (21, 2)
\edge from (22, 2) to (21.8, 2)
\edge from (21.8, 2) to (21.2, 1)
\edge from (21.2, 1) to (21, 1)
\edge from (22, 3) to (21, 3)
\edge from (22, 4) to (21, 4)
\edge from (22, 5) to (21, 5)
\edge from (22, 6) to (21, 6)
\edge from (21, 1) to (20, 1)
\edge from (21, 2) to (20, 2)
\edge from (21, 3) to (20.8, 3)
\edge from (20.8, 3) to (20.2, 4)
\edge from (20.2, 4) to (20, 4)
\edge from (21, 4) to (20.8, 4)
\edge from (20.8, 4) to (20.2, 3)
\edge from (20.2, 3) to (20, 3)
\edge from (21, 5) to (20, 5)
\edge from (21, 6) to (20, 6)
\edge from (20, 1) to (19, 1)
\edge from (20, 2) to (19.8, 2)
\edge from (19.8, 2) to (19.2, 3)
\edge from (19.2, 3) to (19, 3)
\edge from (20, 3) to (19.8, 3)
\edge from (19.8, 3) to (19.2, 2)
\edge from (19.2, 2) to (19, 2)
\edge from (20, 4) to (19, 4)
\edge from (20, 5) to (19, 5)
\edge from (20, 6) to (19, 6)
\edge from (19, 1) to (18, 1)
\edge from (19, 2) to (18, 2)
\edge from (19, 3) to (18, 3)
\edge from (19, 4) to (18, 4)
\edge from (19, 5) to (18.8, 5)
\edge from (18.8, 5) to (18.2, 6)
\edge from (18.2, 6) to (18, 6)
\edge from (19, 6) to (18.8, 6)
\edge from (18.8, 6) to (18.2, 5)
\edge from (18.2, 5) to (18, 5)
\edge from (18, 1) to (17, 1)
\edge from (18, 2) to (17, 2)
\edge from (18, 3) to (17, 3)
\edge from (18, 4) to (17.8, 4)
\edge from (17.8, 4) to (17.2, 5)
\edge from (17.2, 5) to (17, 5)
\edge from (18, 5) to (17.8, 5)
\edge from (17.8, 5) to (17.2, 4)
\edge from (17.2, 4) to (17, 4)
\edge from (18, 6) to (17, 6)
\edge from (17, 1) to (16, 1)
\edge from (17, 2) to (16, 2)
\edge from (17, 3) to (16.8, 3)
\edge from (16.8, 3) to (16.2, 4)
\edge from (16.2, 4) to (16, 4)
\edge from (17, 4) to (16.8, 4)
\edge from (16.8, 4) to (16.2, 3)
\edge from (16.2, 3) to (16, 3)
\edge from (17, 5) to (16, 5)
\edge from (17, 6) to (16, 6)
\edge from (16, 1) to (15, 1)
\edge from (16, 2) to (15.8, 2)
\edge from (15.8, 2) to (15.2, 3)
\edge from (15.2, 3) to (15, 3)
\edge from (16, 3) to (15.8, 3)
\edge from (15.8, 3) to (15.2, 2)
\edge from (15.2, 2) to (15, 2)
\edge from (16, 4) to (15, 4)
\edge from (16, 5) to (15, 5)
\edge from (16, 6) to (15, 6)
\edge from (15, 1) to (14.8, 1)
\edge from (14.8, 1) to (14.2, 2)
\edge from (14.2, 2) to (14, 2)
\edge from (15, 2) to (14.8, 2)
\edge from (14.8, 2) to (14.2, 1)
\edge from (14.2, 1) to (14, 1)
\edge from (15, 3) to (14, 3)
\edge from (15, 4) to (14, 4)
\edge from (15, 5) to (14, 5)
\edge from (15, 6) to (14, 6)
\label {$1$} at (11.5, 1)
\label {$2$} at (11.5, 2)
\label {$3$} at (11.5, 3)
\label {$4$} at (11.5, 4)
\label {$5$} at (11.5, 5)
\label {$6$} at (11.5, 6)
\edge from (11, 1) to (10.8, 1)
\edge from (10.8, 1) to (10.2, 3)
\edge from (10.2, 3) to (10, 3)
\edge from (11, 2) to (10.8, 2)
\edge from (10.8, 2) to (10.2, 2)
\edge from (10.2, 2) to (10, 2)
\edge from (11, 3) to (10.8, 3)
\edge from (10.8, 3) to (10.2, 1)
\edge from (10.2, 1) to (10, 1)
\edge from (11, 4) to (10, 4)
\edge from (11, 5) to (10, 5)
\edge from (11, 6) to (10, 6)
\edge from (10, 1) to (9, 1)
\edge from (10, 2) to (9, 2)
\edge from (10, 3) to (9, 3)
\edge from (10, 4) to (9.8, 4)
\edge from (9.8, 4) to (9.2, 6)
\edge from (9.2, 6) to (9, 6)
\edge from (10, 5) to (9.8, 5)
\edge from (9.8, 5) to (9.2, 5)
\edge from (9.2, 5) to (9, 5)
\edge from (10, 6) to (9.8, 6)
\edge from (9.8, 6) to (9.2, 4)
\edge from (9.2, 4) to (9, 4)
\edge from (9, 1) to (8, 1)
\edge from (9, 2) to (8, 2)
\edge from (9, 3) to (8.8, 3)
\edge from (8.8, 3) to (8.2, 4)
\edge from (8.2, 4) to (8, 4)
\edge from (9, 4) to (8.8, 4)
\edge from (8.8, 4) to (8.2, 3)
\edge from (8.2, 3) to (8, 3)
\edge from (9, 5) to (8, 5)
\edge from (9, 6) to (8, 6)
\edge from (8, 1) to (7, 1)
\edge from (8, 2) to (7.8, 2)
\edge from (7.8, 2) to (7.2, 3)
\edge from (7.2, 3) to (7, 3)
\edge from (8, 3) to (7.8, 3)
\edge from (7.8, 3) to (7.2, 2)
\edge from (7.2, 2) to (7, 2)
\edge from (8, 4) to (7, 4)
\edge from (8, 5) to (7, 5)
\edge from (8, 6) to (7, 6)
\edge from (7, 1) to (6.8, 1)
\edge from (6.8, 1) to (6.2, 2)
\edge from (6.2, 2) to (6, 2)
\edge from (7, 2) to (6.8, 2)
\edge from (6.8, 2) to (6.2, 1)
\edge from (6.2, 1) to (6, 1)
\edge from (7, 3) to (6, 3)
\edge from (7, 4) to (6, 4)
\edge from (7, 5) to (6, 5)
\edge from (7, 6) to (6, 6)
\edge from (6, 1) to (5, 1)
\edge from (6, 2) to (5, 2)
\edge from (6, 3) to (5, 3)
\edge from (6, 4) to (5.8, 4)
\edge from (5.8, 4) to (5.2, 5)
\edge from (5.2, 5) to (5, 5)
\edge from (6, 5) to (5.8, 5)
\edge from (5.8, 5) to (5.2, 4)
\edge from (5.2, 4) to (5, 4)
\edge from (6, 6) to (5, 6)
\edge from (5, 1) to (4, 1)
\edge from (5, 2) to (4, 2)
\edge from (5, 3) to (4.8, 3)
\edge from (4.8, 3) to (4.2, 4)
\edge from (4.2, 4) to (4, 4)
\edge from (5, 4) to (4.8, 4)
\edge from (4.8, 4) to (4.2, 3)
\edge from (4.2, 3) to (4, 3)
\edge from (5, 5) to (4, 5)
\edge from (5, 6) to (4, 6)
\edge from (4, 1) to (3, 1)
\edge from (4, 2) to (3.8, 2)
\edge from (3.8, 2) to (3.2, 3)
\edge from (3.2, 3) to (3, 3)
\edge from (4, 3) to (3.8, 3)
\edge from (3.8, 3) to (3.2, 2)
\edge from (3.2, 2) to (3, 2)
\edge from (4, 4) to (3, 4)
\edge from (4, 5) to (3, 5)
\edge from (4, 6) to (3, 6)
\edge from (3, 1) to (2, 1)
\edge from (3, 2) to (2, 2)
\edge from (3, 3) to (2, 3)
\edge from (3, 4) to (2, 4)
\edge from (3, 5) to (2.8, 5)
\edge from (2.8, 5) to (2.2, 6)
\edge from (2.2, 6) to (2, 6)
\edge from (3, 6) to (2.8, 6)
\edge from (2.8, 6) to (2.2, 5)
\edge from (2.2, 5) to (2, 5)
\edge from (2, 1) to (1, 1)
\edge from (2, 2) to (1, 2)
\edge from (2, 3) to (1, 3)
\edge from (2, 4) to (1.8, 4)
\edge from (1.8, 4) to (1.2, 5)
\edge from (1.2, 5) to (1, 5)
\edge from (2, 5) to (1.8, 5)
\edge from (1.8, 5) to (1.2, 4)
\edge from (1.2, 4) to (1, 4)
\edge from (2, 6) to (1, 6)
\edge from (1, 1) to (0, 1)
\edge from (1, 2) to (0, 2)
\edge from (1, 3) to (0.8, 3)
\edge from (0.8, 3) to (0.2, 4)
\edge from (0.2, 4) to (0, 4)
\edge from (1, 4) to (0.8, 4)
\edge from (0.8, 4) to (0.2, 3)
\edge from (0.2, 3) to (0, 3)
\edge from (1, 5) to (0, 5)
\edge from (1, 6) to (0, 6)
\endgraph}

\nopagebreak
\centerline{Figure 5.  Wiring diagrams for $\A$ (left) and $\A'$
(right)}
\smallskip
Applying the methods described in the previous sections,
we obtain the following braid monodromy generators
$\vec\l=\{\l_k\}$ and $\vec\l'=\{\l_k'\}$:
$$\align
\vec\l&=\{ A_{1,2,3},\ A_{4,5,6},\ A_{1,6}^{\{4,5\}},
\ A_{2,6}^{\{4,5\}},\ A_{3,6}^{\{4,5\}},\ A_{1,5}^{\{4\}},
\ A_{2,5}^{\{4\}},\ A_{3,5}^{\{4\}},\ A_{1,4},\ A_{2,4},\ A_{3,4}\},\\
\vec\l'&=\{ A_{1,2,3},\ A_{1,4,5},\ A_{2,5}^{\{4\}},\ A_{3,5}^{\{4\}},
\ A_{2,4},\ A_{3,4},\ A_{1,6},\ A_{2,6},\ A_{3,6},\ A_{4,6},\ A_{5,6}\}.\\
\endalign$$

Using Proposition~6.9 and some elementary simplifications,
we obtain the following presentations for the groups $G(\A)$ and $G(\A')$:
$$\split G(\A)&=\langle u_1,\dots,u_6 \mid [u_1,u_2,u_3], [u_4,u_5,u_6], [u_1,u_6], [u_2,u_6], [u_3,u_6], [u_1,u_5],\\
&\qquad\qquad\qquad\qquad [u_2,u_5],[u_3,u_5],[u_1,u_4],
[u_2,u_4], [u_3,u_4]\rangle,\\
G(\A')&=\langle v_1,\dots,v_6 \mid [v_1,v_2,v_3], [v_1,v_4,v_5], [v_2,v_5], [v_3,v_5], [v_2,v_4], [v_3,v_4],\\
&\qquad\qquad\qquad\qquad [v_1,v_6],[v_2,v_6],[v_3,v_6],
[v_4,v_6], [v_5,v_6]\rangle,\\
\endsplit$$
These groups are isomorphic.
In fact, one can check that the map $G(\A')\to G(\A)$ defined by
$v_1\mapsto u_1u_5^{-1}u_4^{-1}$, $v_6\mapsto u_4u_5u_6$, $v_i\mapsto u_i$ if
$i\neq 1,6$ is an isomorphism through Tietze-I moves, so the complements 
of $\A$ and $\A'$ are homotopy equivalent.  Falk
obtained analogous results for the original central 3-arrangements by
working with decones as opposed to generic sections.

On the other hand, the lattices of $\A$ and $\A'$ are not isomorphic:
For $\A'$, the two triple points are incident on a line,
while for $\A$, they are not.  By Theorem~7.2,
the monodromies $\l,\l':F_{11}\to B_6$ are not braid-equivalent.
Moreover, the fact that $L(\A) \not\cong L(\A')$ implies, by results
of Jiang and Yau \cite{JY}, that the complements of
$\A$ and $\A'$ are not homeomorphic.

\subhead 7.5   Falk-Sturmfels arrangements
\endsubhead
Consider the pair of (central) plane arrangements
in $\C^3$, with defining polynomials
$$Q^{\pm} = xyz(x+y+z)(x+\gamma y)(y+z)(x+(\gamma +1)y+z)
(-x+\gamma z)(x+\gamma y+\gamma^2z),$$
where $\gamma = (-1 \pm \sqrt{5})/2$.
These arrangements, studied by Falk and Sturmfels (unpublished),
are real realizations of a minimal matroid,
whose realization space is disconnected (see also \cite{BLSWZ}).
They have isomorphic lattices and homotopy equivalent complements.
In fact, Keaty (also unpublished) has shown that the oriented
matroids of these arrangements are isomorphic.  Thus, by results of
Bj\"orner and Ziegler \cite{BZ}, their complements are homeomorphic.
Checking that $H^+=\{z=1-{2\over 5}x+{2\over 7}y\}$ and
$H^-=\{z=1-{4\over 9}x+{1\over 7}y\}$ are generic with respect
to these arrangements, we get a pair of real line arrangements,
$\A^\pm$, in $\C^2$ by taking sections.  Wiring diagrams for
these line arrangements are depicted below.

\bigskip

\centerline{
\beginsmallgraph
\label {$1$} at (27.5, 1)
\label {$2$} at (27.5, 2)
\label {$3$} at (27.5, 3)
\label {$4$} at (27.5, 4)
\label {$5$} at (27.5, 5)
\label {$6$} at (27.5, 6)
\label {$7$} at (27.5, 7)
\label {$8$} at (27.5, 8)
\label {$9$} at (27.5, 9)
\edge from (27, 1) to (26, 1)
\edge from (27, 2) to (26, 2)
\edge from (27, 3) to (26, 3)
\edge from (27, 4) to (26, 4)
\edge from (27, 5) to (26, 5)
\edge from (27, 6) to (26.8, 6)
\edge from (26.8, 6) to (26.2, 7)
\edge from (26.2, 7) to (26, 7)
\edge from (27, 7) to (26.8, 7)
\edge from (26.8, 7) to (26.2, 6)
\edge from (26.2, 6) to (26, 6)
\edge from (27, 8) to (26, 8)
\edge from (27, 9) to (26, 9)
\edge from (26, 1) to (25, 1)
\edge from (26, 2) to (25.8, 2)
\edge from (25.8, 2) to (25.2, 3)
\edge from (25.2, 3) to (25, 3)
\edge from (26, 3) to (25.8, 3)
\edge from (25.8, 3) to (25.2, 2)
\edge from (25.2, 2) to (25, 2)
\edge from (26, 4) to (25, 4)
\edge from (26, 5) to (25.8, 5)
\edge from (25.8, 5) to (25.2, 6)
\edge from (25.2, 6) to (25, 6)
\edge from (26, 6) to (25.8, 6)
\edge from (25.8, 6) to (25.2, 5)
\edge from (25.2, 5) to (25, 5)
\edge from (26, 7) to (25, 7)
\edge from (26, 8) to (25, 8)
\edge from (26, 9) to (25, 9)
\edge from (25, 1) to (24, 1)
\edge from (25, 2) to (24, 2)
\edge from (25, 3) to (24.8, 3)
\edge from (24.8, 3) to (24.2, 5)
\edge from (24.2, 5) to (24, 5)
\edge from (25, 4) to (24.8, 4)
\edge from (24.8, 4) to (24.2, 4)
\edge from (24.2, 4) to (24, 4)
\edge from (25, 5) to (24.8, 5)
\edge from (24.8, 5) to (24.2, 3)
\edge from (24.2, 3) to (24, 3)
\edge from (25, 6) to (24, 6)
\edge from (25, 7) to (24, 7)
\edge from (25, 8) to (24, 8)
\edge from (25, 9) to (24, 9)
\edge from (24, 1) to (23, 1)
\edge from (24, 2) to (23, 2)
\edge from (24, 3) to (23, 3)
\edge from (24, 4) to (23, 4)
\edge from (24, 5) to (23.8, 5)
\edge from (23.8, 5) to (23.2, 7)
\edge from (23.2, 7) to (23, 7)
\edge from (24, 6) to (23.8, 6)
\edge from (23.8, 6) to (23.2, 6)
\edge from (23.2, 6) to (23, 6)
\edge from (24, 7) to (23.8, 7)
\edge from (23.8, 7) to (23.2, 5)
\edge from (23.2, 5) to (23, 5)
\edge from (24, 8) to (23, 8)
\edge from (24, 9) to (23, 9)
\edge from (23, 1) to (22.8, 1)
\edge from (22.8, 1) to (22.2, 3)
\edge from (22.2, 3) to (22, 3)
\edge from (23, 2) to (22, 2)
\edge from (23, 3) to (22.8, 3)
\edge from (22.8, 3) to (22.2, 1)
\edge from (22.2, 1) to (22, 1)
\edge from (23, 2) to (22, 2)
\edge from (23, 4) to (22, 4)
\edge from (23, 5) to (22, 5)
\edge from (23, 6) to (22, 6)
\edge from (23, 7) to (22.8, 7)
\edge from (22.8, 7) to (22.2, 8)
\edge from (22.2, 8) to (22, 8)
\edge from (23, 8) to (22.8, 8)
\edge from (22.8, 8) to (22.2, 7)
\edge from (22.2, 7) to (22, 7)
\edge from (23, 9) to (22, 9)
\edge from (22, 1) to (21, 1)
\edge from (22, 2) to (21, 2)
\edge from (22, 3) to (21.8, 3)
\edge from (21.8, 3) to (21.2, 5)
\edge from (21.2, 5) to (21, 5)
\edge from (22, 4) to (21.8, 4)
\edge from (21.8, 4) to (21.2, 4)
\edge from (21.2, 4) to (21, 4)
\edge from (22, 5) to (21.8, 5)
\edge from (21.8, 5) to (21.2, 3)
\edge from (21.2, 3) to (21, 3)
\edge from (22, 6) to (21, 6)
\edge from (22, 7) to (21, 7)
\edge from (22, 8) to (21, 8)
\edge from (22, 9) to (21, 9)
\edge from (21, 1) to (20, 1)
\edge from (21, 2) to (20, 2)
\edge from (21, 3) to (20, 3)
\edge from (21, 4) to (20, 4)
\edge from (21, 5) to (20.8, 5)
\edge from (20.8, 5) to (20.2, 7)
\edge from (20.2, 7) to (20, 7)
\edge from (21, 6) to (20.8, 6)
\edge from (20.8, 6) to (20.2, 6)
\edge from (20.2, 6) to (20, 6)
\edge from (21, 7) to (20.8, 7)
\edge from (20.8, 7) to (20.2, 5)
\edge from (20.2, 5) to (20, 5)
\edge from (21, 8) to (20, 8)
\edge from (21, 9) to (20, 9)
\edge from (20, 1) to (19, 1)
\edge from (20, 2) to (19, 2)
\edge from (20, 3) to (19, 3)
\edge from (20, 4) to (19.8, 4)
\edge from (19.8, 4) to (19.2, 5)
\edge from (19.2, 5) to (19, 5)
\edge from (20, 5) to (19.8, 5)
\edge from (19.8, 5) to (19.2, 4)
\edge from (19.2, 4) to (19, 4)
\edge from (20, 6) to (19, 6)
\edge from (20, 7) to (19, 7)
\edge from (20, 8) to (19, 8)
\edge from (20, 9) to (19, 9)
\edge from (19, 1) to (18, 1)
\edge from (19, 2) to (18.8, 2)
\edge from (18.8, 2) to (18.2, 4)
\edge from (18.2, 4) to (18, 4)
\edge from (19, 3) to (18, 3)
\edge from (19, 4) to (18.8, 4)
\edge from (18.8, 4) to (18.2, 2)
\edge from (18.2, 2) to (18, 2)
\edge from (19, 5) to (18, 5)
\edge from (19, 6) to (18, 6)
\edge from (19, 7) to (18.8, 7)
\edge from (18.8, 7) to (18.2, 9)
\edge from (18.2, 9) to (18, 9)
\edge from (19, 8) to (18.8, 8)
\edge from (18.8, 8) to (18.2, 8)
\edge from (18.2, 8) to (18, 8)
\edge from (19, 9) to (18.8, 9)
\edge from (18.8, 9) to (18.2, 7)
\edge from (18.2, 7) to (18, 7)
\edge from (18, 1) to (17, 1)
\edge from (18, 2) to (17, 2)
\edge from (18, 3) to (17, 3)
\edge from (18, 4) to (17.8, 4)
\edge from (17.8, 4) to (17.2, 7)
\edge from (17.2, 7) to (17, 7)
\edge from (18, 5) to (17.8, 5)
\edge from (17.8, 5) to (17.2, 6)
\edge from (17.2, 6) to (17, 6)
\edge from (18, 6) to (17.8, 6)
\edge from (17.8, 6) to (17.2, 5)
\edge from (17.2, 5) to (17, 5)
\edge from (18, 7) to (17.8, 7)
\edge from (17.8, 7) to (17.2, 4)
\edge from (17.2, 4) to (17, 4)
\edge from (18, 8) to (17, 8)
\edge from (18, 9) to (17, 9)
\edge from (17, 1) to (16, 1)
\edge from (17, 2) to (16, 2)
\edge from (17, 3) to (16.8, 3)
\edge from (16.8, 3) to (16.2, 4)
\edge from (16.2, 4) to (16, 4)
\edge from (17, 4) to (16.8, 4)
\edge from (16.8, 4) to (16.2, 3)
\edge from (16.2, 3) to (16, 3)
\edge from (17, 5) to (16, 5)
\edge from (17, 6) to (16, 6)
\edge from (17, 7) to (16, 7)
\edge from (17, 8) to (16, 8)
\edge from (17, 9) to (16, 9)
\edge from (16, 1) to (15.8, 1)
\edge from (15.8, 1) to (15.2, 3)
\edge from (15.2, 3) to (15, 3)
\edge from (16, 2) to (15.8, 2)
\edge from (15.8, 2) to (15.2, 2)
\edge from (15.2, 2) to (15, 2)
\edge from (16, 3) to (15.8, 3)
\edge from (15.8, 3) to (15.2, 1)
\edge from (15.2, 1) to (15, 1)
\edge from (16, 4) to (15, 4)
\edge from (16, 5) to (15, 5)
\edge from (16, 6) to (15, 6)
\edge from (16, 7) to (15, 7)
\edge from (16, 8) to (15, 8)
\edge from (16, 9) to (15, 9)
\label {$1$} at (13.5, 1)
\label {$2$} at (13.5, 2)
\label {$3$} at (13.5, 3)
\label {$4$} at (13.5, 4)
\label {$5$} at (13.5, 5)
\label {$6$} at (13.5, 6)
\label {$7$} at (13.5, 7)
\label {$8$} at (13.5, 8)
\label {$9$} at (13.5, 9)
\edge from (13, 1) to (12.8, 1)
\edge from (12.8, 1) to (12.2, 3)
\edge from (12.2, 3) to (12, 3)
\edge from (13, 2) to (12.8, 2)
\edge from (12.8, 2) to (12.2, 2)
\edge from (12.2, 2) to (12, 2)
\edge from (13, 3) to (12.8, 3)
\edge from (12.8, 3) to (12.2, 1)
\edge from (12.2, 1) to (12, 1)
\edge from (13, 4) to (12, 4)
\edge from (13, 5) to (12, 5)
\edge from (13, 6) to (12, 6)
\edge from (13, 7) to (12, 7)
\edge from (13, 8) to (12, 8)
\edge from (13, 9) to (12, 9)
\edge from (12, 1) to (11, 1)
\edge from (12, 2) to (11, 2)
\edge from (12, 3) to (11.8, 3)
\edge from (11.8, 3) to (11.2, 4)
\edge from (11.2, 4) to (11, 4)
\edge from (12, 4) to (11.8, 4)
\edge from (11.8, 4) to (11.2, 3)
\edge from (11.2, 3) to (11, 3)
\edge from (12, 5) to (11, 5)
\edge from (12, 6) to (11, 6)
\edge from (12, 7) to (11, 7)
\edge from (12, 8) to (11, 8)
\edge from (12, 9) to (11, 9)
\edge from (11, 1) to (10, 1)
\edge from (11, 2) to (10, 2)
\edge from (11, 3) to (10, 3)
\edge from (11, 4) to (10.8, 4)
\edge from (10.8, 4) to (10.2, 6)
\edge from (10.2, 6) to (10, 6)
\edge from (11, 5) to (10.8, 5)
\edge from (10.8, 5) to (10.2, 5)
\edge from (10.2, 5) to (10, 5)
\edge from (11, 6) to (10.8, 6)
\edge from (10.8, 6) to (10.2, 4)
\edge from (10.2, 4) to (10, 4)
\edge from (11, 7) to (10, 7)
\edge from (11, 8) to (10, 8)
\edge from (11, 9) to (10, 9)
\edge from (10, 1) to (9, 1)
\edge from (10, 2) to (9, 2)
\edge from (10, 3) to (9, 3)
\edge from (10, 4) to (9, 4)
\edge from (10, 5) to (9, 5)
\edge from (10, 6) to (9.8, 6)
\edge from (9.8, 6) to (9.2, 7)
\edge from (9.2, 7) to (9, 7)
\edge from (10, 7) to (9.8, 7)
\edge from (9.8, 7) to (9.2, 6)
\edge from (9.2, 6) to (9, 6)
\edge from (10, 8) to (9, 8)
\edge from (10, 9) to (9, 9)
\edge from (9, 1) to (8, 1)
\edge from (9, 2) to (8.8, 2)
\edge from (8.8, 2) to (8.2, 4)
\edge from (8.2, 4) to (8, 4)
\edge from (9, 3) to (8, 3)
\edge from (9, 4) to (8.8, 4)
\edge from (8.8, 4) to (8.2, 2)
\edge from (8.2, 2) to (8, 2)
\edge from (9, 5) to (8, 5)
\edge from (9, 6) to (8, 6)
\edge from (9, 7) to (8.8, 7)
\edge from (8.8, 7) to (8.2, 9)
\edge from (8.2, 9) to (8, 9)
\edge from (9, 8) to (8.8, 8)
\edge from (8.8, 8) to (8.2, 8)
\edge from (8.2, 8) to (8, 8)
\edge from (9, 9) to (8.8, 9)
\edge from (8.8, 9) to (8.2, 7)
\edge from (8.2, 7) to (8, 7)
\edge from (8, 1) to (7, 1)
\edge from (8, 2) to (7, 2)
\edge from (8, 3) to (7, 3)
\edge from (8, 4) to (7.8, 4)
\edge from (7.8, 4) to (7.2, 7)
\edge from (7.2, 7) to (7, 7)
\edge from (8, 5) to (7.8, 5)
\edge from (7.8, 5) to (7.2, 6)
\edge from (7.2, 6) to (7, 6)
\edge from (8, 6) to (7.8, 6)
\edge from (7.8, 6) to (7.2, 5)
\edge from (7.2, 5) to (7, 5)
\edge from (8, 7) to (7.8, 7)
\edge from (7.8, 7) to (7.2, 4)
\edge from (7.2, 4) to (7, 4)
\edge from (8, 8) to (7, 8)
\edge from (8, 9) to (7, 9)
\edge from (7, 1) to (6, 1)
\edge from (7, 2) to (6, 2)
\edge from (7, 3) to (6.8, 3)
\edge from (6.8, 3) to (6.2, 4)
\edge from (6.2, 4) to (6, 4)
\edge from (7, 4) to (6.8, 4)
\edge from (6.8, 4) to (6.2, 3)
\edge from (6.2, 3) to (6, 3)
\edge from (7, 5) to (6, 5)
\edge from (7, 6) to (6, 6)
\edge from (7, 7) to (6, 7)
\edge from (7, 8) to (6, 8)
\edge from (7, 9) to (6, 9)
\edge from (6, 1) to (5.8, 1)
\edge from (5.8, 1) to (5.2, 3)
\edge from (5.2, 3) to (5, 3)
\edge from (6, 2) to (5.8, 2)
\edge from (5.8, 2) to (5.2, 2)
\edge from (5.2, 2) to (5, 2)
\edge from (6, 3) to (5.8, 3)
\edge from (5.8, 3) to (5.2, 1)
\edge from (5.2, 1) to (5, 1)
\edge from (6, 4) to (5, 4)
\edge from (6, 5) to (5, 5)
\edge from (6, 6) to (5, 6)
\edge from (6, 7) to (5, 7)
\edge from (6, 8) to (5, 8)
\edge from (6, 9) to (5, 9)
\edge from (5, 1) to (4, 1)
\edge from (5, 2) to (4, 2)
\edge from (5, 3) to (4.8, 3)
\edge from (4.8, 3) to (4.2, 5)
\edge from (4.2, 5) to (4, 5)
\edge from (5, 4) to (4.8, 4)
\edge from (4.8, 4) to (4.2, 4)
\edge from (4.2, 4) to (4, 4)
\edge from (5, 5) to (4.8, 5)
\edge from (4.8, 5) to (4.2, 3)
\edge from (4.2, 3) to (4, 3)
\edge from (5, 6) to (4, 6)
\edge from (5, 7) to (4, 7)
\edge from (5, 8) to (4, 8)
\edge from (5, 9) to (4, 9)
\edge from (4, 1) to (3, 1)
\edge from (4, 2) to (3, 2)
\edge from (4, 3) to (3, 3)
\edge from (4, 4) to (3, 4)
\edge from (4, 5) to (3.8, 5)
\edge from (3.8, 5) to (3.2, 6)
\edge from (3.2, 6) to (3, 6)
\edge from (4, 6) to (3.8, 6)
\edge from (3.8, 6) to (3.2, 5)
\edge from (3.2, 5) to (3, 5)
\edge from (4, 7) to (3.8, 7)
\edge from (3.8, 7) to (3.2, 8)
\edge from (3.2, 8) to (3, 8)
\edge from (4, 8) to (3.8, 8)
\edge from (3.8, 8) to (3.2, 7)
\edge from (3.2, 7) to (3, 7)
\edge from (4, 9) to (3, 9)
\edge from (3, 1) to (2, 1)
\edge from (3, 2) to (2, 2)
\edge from (3, 3) to (2, 3)
\edge from (3, 4) to (2, 4)
\edge from (3, 5) to (2, 5)
\edge from (3, 6) to (2.8, 6)
\edge from (2.8, 6) to (2.2, 7)
\edge from (2.2, 7) to (2, 7)
\edge from (3, 7) to (2.8, 7)
\edge from (2.8, 7) to (2.2, 6)
\edge from (2.2, 6) to (2, 6)
\edge from (3, 8) to (2, 8)
\edge from (3, 9) to (2, 9)
\edge from (2, 1) to (1, 1)
\edge from (2, 2) to (1, 2)
\edge from (2, 3) to (1, 3)
\edge from (2, 4) to (1.8, 4)
\edge from (1.8, 4) to (1.2, 6)
\edge from (1.2, 6) to (1, 6)
\edge from (2, 5) to (1.8, 5)
\edge from (1.8, 5) to (1.2, 5)
\edge from (1.2, 5) to (1, 5)
\edge from (2, 6) to (1.8, 6)
\edge from (1.8, 6) to (1.2, 4)
\edge from (1.2, 4) to (1, 4)
\edge from (2, 7) to (1, 7)
\edge from (2, 8) to (1, 8)
\edge from (2, 9) to (1, 9)
\edge from (1, 1) to (0, 1)
\edge from (1, 2) to (0.8, 2)
\edge from (0.8, 2) to (0.2, 4)
\edge from (0.2, 4) to (0, 4)
\edge from (1, 3) to (0.8, 3)
\edge from (0.8, 3) to (0.2, 3)
\edge from (0.2, 3) to (0, 3)
\edge from (1, 4) to (0.8, 4)
\edge from (0.8, 4) to (0.2, 2)
\edge from (0.2, 2) to (0, 2)
\edge from (1, 5) to (0, 5)
\edge from (1, 6) to (0, 6)
\edge from (1, 7) to (0, 7)
\edge from (1, 8) to (0, 8)
\edge from (1, 9) to (0, 9)
\endgraph}

\nopagebreak
\centerline{Figure 6.  Wiring diagrams for $\A^+$ (left)
and $\A^-$ (right)}
\smallskip
Applying the techniques described in the previous sections,
we obtain the following braid monodromy generators:
$$\align
\vec\l^+&=\{ A_{1,2,3},\ A_{1,4},\ A_{1,5,6},\ A_{1,7},\
A_{1,8,9},\
A_{2,4,6}^{\{5\}},\ A_{2,5,7,9}^{\{8\}},\ A_{4,9}^{\{5,7,8\}},\
A_{3,6,9}^{\{4,5,7,8\}},\ A_{3,4,7}^{\{5\}},\\ &\qquad\qquad\qquad
A_{2,8},\ A_{3,5},\ A_{3,8},\ A_{4,5,8},\ A_{6,7,8}
 \},\\
\vec\l^-&=\{ A_{6,7},\ A_{5,7}^{\{6\}}\ ,A_{2,3},\
A_{2,4,7}^{\{5,6\}},\ A_{2,5,6},\
A_{2,8},\ A_{1,3,7}^{\{2,4,5,6\}},\
A_{1,4,6}^{\{2,5\}},\ A_{1,5,8}^{\{2\}},\
A_{4,8}^{\{5\}},\ A_{1,2,9},\\ &\qquad\qquad\qquad
A_{3,6,8}^{\{4,5\}},\ A_{3,4,5,9},\ A_{6,9},\ A_{7,8,9}
\}.\\
\endalign$$

The monodromy homomorphisms $\l^+, \l^-:F_{15}\to B_9$ are braid-equivalent.
If $I=\{i,i+1,\dots,j\}$, write $\mu_{i,j}=\mu_I$, see \thetag{8}.
One can check that $\l^+\circ\psi=\conj_{\phi}\circ\l^-$, where
$\psi\in B_{15}$ is given by
$$\psi=(\s_1\s_2\cdots\s_{14})^4\s_6
\mu_{3,6}
\s_6\s_{10}\s_{11}\s_{12}\s_{13}
\mu_{7,10}
\s_{10}
\mu_{11,15}
\s_4\s_5\s_6\s_9\s_{10}
\mu_{6,9}
\s_5\s_{11}$$
and $\phi\in B_9$ is given by
$\phi = (\s_8\s_7 \cdots\s_1)^4\s_3\s_4\s_3\s_2\s_5\s_6$.
It follows from Remark 4.3 that the groups $G(\A^+)$ and $G(\A^-)$ are
isomorphic.

\example{Remark 7.6}  If $\lambda:F_s\to P_n$ is the braid monodromy
of an arrangement or wiring diagram, let $\Gamma=\im(\lambda) < P_n$.
In \cite{CF}, it is asserted that if $\W$ and $\W'$ are wiring diagrams
determining the {\it same} underlying matroid, then the corresponding
braid monodromy subgroups $\Gamma$ and $\Gamma'$ of $P_n$ are equal.
A subsequent result for real arrangements may be found in \cite{Cor}.

These results cannot be strengthened.  By construction, the wiring diagrams
$\W^+$ and $\W^-$ above determine {\it isomorphic}
underlying matroids.  However, their braid monodromy groups
$\Gamma^+$ and $\Gamma^-$ do not coincide.
In fact, $\Gamma^+$ and $\Gamma^-$ are not conjugate in $P_9$
(although, as shown previously, they are conjugate in $B_9$).  This
can be seen by using the representation $\theta: P_9 \to \GL(8,\Z\Z^{10})$,
which is obtained from the generalized Gassner representation
$\hat\theta^1_{9,2,2}$
of \cite{CS1},~Section 5.8, by restriction to a direct summand.
The corresponding modules of coinvariants, $A_{\theta}(\A^\pm)=H_0(F_{15};
(\Z\Z^{10})^8_{\theta\circ\lambda^{\pm}})$, are not isomorphic.
One can show that the graded modules associated to their $I$-adic completions
have different Hilbert series.

This difference may be explained combinatorially as follows.  Though
the (little) oriented matroids of $\W^+$ and $\W^-$ are isomorphic,
their big oriented matroids are not.  One can check that the spectra
of the tope graphs (see \cite{BSLWZ}) associated to these big oriented
matroids differ.
\endexample

\subhead 7.7 MacLane configurations
\endsubhead
We consider complex conjugate arrangements arising from the MacLane
matroid $\text{ML}_8$ \cite{MacL}.  This matroid is minimal
non-orientable in the sense that it is the smallest matroid that is
realizable over $\C$ but not over $\R$ \cite{BLSWZ}.
Following Rybnikov \cite{Ry}, we take arrangements with defining
polynomials
$$Q^\pm = xyz(y-x)(z-x)(z+\omega y)
(z+\omega^2 x+\omega y)(z-x-\omega^2 y),$$
where $\omega=(-1\pm \sqrt{-3})/2$, as complex realizations of
this matroid.  Deconing by setting $x=1$, we obtain two affine
arrangements $\A^\pm$ in $\C^2$.  By construction, $\A^+$ and $\A^-$
are lattice-isomorphic.   Also, note that $\A^+$ and $\A^-$ are
conjugate arrangements; in particular, they have diffeomorphic complements,
and thus isomorphic groups.

Check that the projection $\pi(y,z) = 3y+z$ is generic with
respect to $\A^+$. Changing coordinates accordingly, and choosing
an admissible path $\xi$,
we obtain the braided wiring diagram $\W^+$ depicted below.

\bigskip

\centerline{
\beginsmallgraph
\label {$1$} at (19.5, 1)
\label {$2$} at (19.5, 2)
\label {$3$} at (19.5, 3)
\label {$4$} at (19.5, 4)
\label {$5$} at (19.5, 5)
\label {$6$} at (19.5, 6)
\label {$7$} at (19.5, 7)
\edge from (19, 1) to (18, 1)
\edge from (19, 2) to (18, 2)
\edge from (19, 3) to (18, 3)
\edge from (19, 4) to (18.8, 4)
\edge from (18.8, 4) to (18.2, 6)
\edge from (18.2, 6) to (18, 6)
\edge from (19, 5) to (18.8, 5)
\edge from (18.8, 5) to (18.2, 5)
\edge from (18.2, 5) to (18, 5)
\edge from (19, 6) to (18.8, 6)
\edge from (18.8, 6) to (18.2, 4)
\edge from (18.2, 4) to (18, 4)
\edge from (19, 7) to (18, 7)
\edge from (18, 1) to (17, 1)
\edge from (18, 2) to (17, 2)
\edge from (18, 3) to (17.8, 3)
\edge from (17.8, 3) to (17.2, 4)
\edge from (17.2, 4) to (17, 4)
\edge from (18, 4) to (17.8, 4)
\edge from (17.8, 4) to (17.2, 3)
\edge from (17.2, 3) to (17, 3)
\edge from (18, 5) to (17, 5)
\edge from (18, 6) to (17, 6)
\edge from (18, 7) to (17, 7)
\edge from (17, 1) to (16.8, 1)
\edge from (16.8, 1) to (16.2, 3)
\edge from (16.2, 3) to (16, 3)
\edge from (17, 2) to (16.8, 2)
\edge from (16.8, 2) to (16.2, 2)
\edge from (16.2, 2) to (16, 2)
\edge from (17, 3) to (16.8, 3)
\edge from (16.8, 3) to (16.2, 1)
\edge from (16.2, 1) to (16, 1)
\edge from (17, 4) to (16, 4)
\edge from (17, 5) to (16, 5)
\edge from (17, 6) to (16, 6)
\edge from (17, 7) to (16, 7)
\edge from (16, 1) to (15, 1)
\edge from (16, 2) to (15, 2)
\edge from (16, 3) to (15, 3)
\edge from (16, 4) to (15, 4)
\edge from (16, 5) to (15.8, 5)
\edge from (15.8, 5) to (15.2, 6)
\edge from (15.2, 6) to (15, 6)
\edge from (16, 6) to (15.8, 6)
\edge from (15.8, 6) to (15.6, 5.67)
\edge from (15.4, 5.33) to (15.2, 5)
\edge from (15.2, 5) to (15, 5)
\edge from (16, 7) to (15, 7)
\edge from (15, 1) to (14, 1)
\edge from (15, 2) to (14, 2)
\edge from (15, 3) to (14, 3)
\edge from (15, 4) to (14.8, 4)
\edge from (14.8, 4) to (14.6, 4.33)
\edge from (14.4, 4.67) to (14.2, 5)
\edge from (14.2, 5) to (14, 5)
\edge from (15, 5) to (14.8, 5)
\edge from (14.8, 5) to (14.2, 4)
\edge from (14.2, 4) to (14, 4)
\edge from (15, 6) to (14, 6)
\edge from (15, 7) to (14, 7)
\edge from (14, 1) to (13, 1)
\edge from (14, 2) to (13, 2)
\edge from (14, 3) to (13.8, 3)
\edge from (13.8, 3) to (13.2, 5)
\edge from (13.2, 5) to (13, 5)
\edge from (14, 4) to (13.8, 4)
\edge from (13.8, 4) to (13.2, 4)
\edge from (13.2, 4) to (13, 4)
\edge from (14, 5) to (13.8, 5)
\edge from (13.8, 5) to (13.2, 3)
\edge from (13.2, 3) to (13, 3)
\edge from (14, 6) to (13, 6)
\edge from (14, 7) to (13, 7)
\edge from (13, 1) to (12, 1)
\edge from (13, 2) to (12, 2)
\edge from (13, 3) to (12, 3)
\edge from (13, 4) to (12.8, 4)
\edge from (12.8, 4) to (12.6, 4.33)
\edge from (12.4, 4.67) to (12.2, 5)
\edge from (12.2, 5) to (12, 5)
\edge from (13, 5) to (12.8, 5)
\edge from (12.8, 5) to (12.2, 4)
\edge from (12.2, 4) to (12, 4)
\edge from (13, 6) to (12, 6)
\edge from (13, 7) to (12, 7)
\edge from (12, 1) to (11, 1)
\edge from (12, 2) to (11, 2)
\edge from (12, 3) to (11.8, 3)
\edge from (11.8, 3) to (11.6, 3.33)
\edge from (11.4, 3.67) to (11.2, 4)
\edge from (11.2, 4) to (11, 4)
\edge from (12, 4) to (11.8, 4)
\edge from (11.8, 4) to (11.2, 3)
\edge from (11.2, 3) to (11, 3)
\edge from (12, 5) to (11, 5)
\edge from (12, 6) to (11, 6)
\edge from (12, 7) to (11, 7)
\edge from (11, 1) to (10, 1)
\edge from (11, 2) to (10, 2)
\edge from (11, 3) to (10, 3)
\edge from (11, 4) to (10, 4)
\edge from (11, 5) to (10.8, 5)
\edge from (10.8, 5) to (10.6, 5.33)
\edge from (10.4, 5.67) to (10.2, 6)
\edge from (10.2, 6) to (10, 6)
\edge from (11, 6) to (10.8, 6)
\edge from (10.8, 6) to (10.2, 5)
\edge from (10.2, 5) to (10, 5)
\edge from (11, 7) to (10, 7)
\edge from (10, 1) to (9, 1)
\edge from (10, 2) to (9.8, 2)
\edge from (9.8, 2) to (9.6, 2.33)
\edge from (9.4, 2.67) to (9.2, 3)
\edge from (9.2, 3) to (9, 3)
\edge from (10, 3) to (9.8, 3)
\edge from (9.8, 3) to (9.2, 2)
\edge from (9.2, 2) to (9, 2)
\edge from (10, 4) to (9, 4)
\edge from (10, 5) to (9, 5)
\edge from (10, 6) to (9, 6)
\edge from (10, 7) to (9, 7)
\edge from (9, 1) to (8, 1)
\edge from (9, 2) to (8, 2)
\edge from (9, 3) to (8.8, 3)
\edge from (8.8, 3) to (8.2, 4)
\edge from (8.2, 4) to (8, 4)
\edge from (9, 4) to (8.8, 4)
\edge from (8.8, 4) to (8.6, 3.67)
\edge from (8.4, 3.33) to (8.2, 3)
\edge from (8.2, 3) to (8, 3)
\edge from (9, 5) to (8, 5)
\edge from (9, 6) to (8, 6)
\edge from (9, 7) to (8, 7)
\edge from (8, 1) to (7, 1)
\edge from (8, 2) to (7, 2)
\edge from (8, 3) to (7, 3)
\edge from (8, 4) to (7.8, 4)
\edge from (7.8, 4) to (7.2, 5)
\edge from (7.2, 5) to (7, 5)
\edge from (8, 5) to (7.8, 5)
\edge from (7.8, 5) to (7.2, 4)
\edge from (7.2, 4) to (7, 4)
\edge from (8, 6) to (7, 6)
\edge from (8, 7) to (7, 7)
\edge from (7, 1) to (6, 1)
\edge from (7, 2) to (6, 2)
\edge from (7, 3) to (6, 3)
\edge from (7, 4) to (6.8, 4)
\edge from (6.8, 4) to (6.2, 5)
\edge from (6.2, 5) to (6, 5)
\edge from (7, 5) to (6.8, 5)
\edge from (6.8, 5) to (6.6, 4.67)
\edge from (6.4, 4.33) to (6.2, 4)
\edge from (6.2, 4) to (6, 4)
\edge from (7, 6) to (6, 6)
\edge from (7, 7) to (6, 7)
\edge from (6, 1) to (5, 1)
\edge from (6, 2) to (5, 2)
\edge from (6, 3) to (5.8, 3)
\edge from (5.8, 3) to (5.6, 3.33)
\edge from (5.4, 3.67) to (5.2, 4)
\edge from (5.2, 4) to (5, 4)
\edge from (6, 4) to (5.8, 4)
\edge from (5.8, 4) to (5.2, 3)
\edge from (5.2, 3) to (5, 3)
\edge from (6, 5) to (5, 5)
\edge from (6, 6) to (5, 6)
\edge from (6, 7) to (5, 7)
\edge from (5, 1) to (4, 1)
\edge from (5, 2) to (4.8, 2)
\edge from (4.8, 2) to (4.2, 3)
\edge from (4.2, 3) to (4, 3)
\edge from (5, 3) to (4.8, 3)
\edge from (4.8, 3) to (4.6, 2.67)
\edge from (4.4, 2.33) to (4.2, 2)
\edge from (4.2, 2) to (4, 2)
\edge from (5, 4) to (4, 4)
\edge from (5, 5) to (4, 5)
\edge from (5, 6) to (4, 6)
\edge from (5, 7) to (4, 7)
\edge from (4, 1) to (3, 1)
\edge from (4, 2) to (3, 2)
\edge from (4, 3) to (3, 3)
\edge from (4, 4) to (3, 4)
\edge from (4, 5) to (3, 5)
\edge from (4, 6) to (3.8, 6)
\edge from (3.8, 6) to (3.2, 7)
\edge from (3.2, 7) to (3, 7)
\edge from (4, 7) to (3.8, 7)
\edge from (3.8, 7) to (3.2, 6)
\edge from (3.2, 6) to (3, 6)
\edge from (3, 1) to (2, 1)
\edge from (3, 2) to (2, 2)
\edge from (3, 3) to (2.8, 3)
\edge from (2.8, 3) to (2.2, 4)
\edge from (2.2, 4) to (2, 4)
\edge from (3, 4) to (2.8, 4)
\edge from (2.8, 4) to (2.6, 3.67)
\edge from (2.4, 3.33) to (2.2, 3)
\edge from (2.2, 3) to (2, 3)
\edge from (3, 5) to (2, 5)
\edge from (3, 6) to (2, 6)
\edge from (3, 7) to (2, 7)
\edge from (2, 1) to (1, 1)
\edge from (2, 2) to (1, 2)
\edge from (2, 3) to (1, 3)
\edge from (2, 4) to (1.8, 4)
\edge from (1.8, 4) to (1.2, 6)
\edge from (1.2, 6) to (1, 6)
\edge from (2, 5) to (1.8, 5)
\edge from (1.8, 5) to (1.2, 5)
\edge from (1.2, 5) to (1, 5)
\edge from (2, 6) to (1.8, 6)
\edge from (1.8, 6) to (1.2, 4)
\edge from (1.2, 4) to (1, 4)
\edge from (2, 7) to (1, 7)
\edge from (1, 1) to (0, 1)
\edge from (1, 2) to (0.8, 2)
\edge from (0.8, 2) to (0.2, 4)
\edge from (0.2, 4) to (0, 4)
\edge from (1, 3) to (0.8, 3)
\edge from (0.8, 3) to (0.2, 3)
\edge from (0.2, 3) to (0, 3)
\edge from (1, 4) to (0.8, 4)
\edge from (0.8, 4) to (0.2, 2)
\edge from (0.2, 2) to (0, 2)
\edge from (1, 5) to (0, 5)
\edge from (1, 6) to (0, 6)
\edge from (1, 7) to (0, 7)
\endgraph
}

\medskip
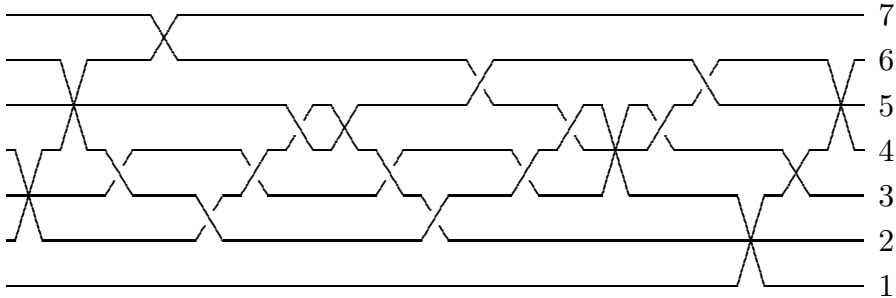
\centerline{Figure 7.  Braided wiring diagram for $\A^+$}
\smallskip

\newpage

\noindent From the braided wiring diagram $\W^+$, we see that
$$\matrix
I_1=\{4,5,6\},\hfill&\b_{1,2} = 1,\hfill&
I_2=\{3,4\},\hfill&\b_{2,3}=1,\hfill\\
I_3=\{1,2,3\},\hfill&\b_{3,4}=\s_4\s_5^{-1},\hfill&
I_4=\{3,4,5\},\hfill&\b_{4,5}=\s_3^{-1}\s_2\s_5\s_3\s_4,\hfill\\
I_5=\{4,5\},\hfill&\b_{5,6}=\s_2^{-1}\s_3\s_4^{-1},\hfill&
I_6=\{6,7\},\hfill&\b_{6,7}=\s_3^{-1},\hfill\\
I_7=\{4,5,6\},\hfill&\b_{7,8}=1,\hfill&I_8=\{2,3,4\}.\hfill
\endmatrix$$
Since $\A^+$ and $\A^-$ are conjugate, a braided wiring diagram
$\W^-$ for $\A^-$ may be obtained from $\W^+$ by switching the
crossings of the intermediate braids, as noted in~5.4.
Applying the algorithm of 5.3 and carrying out some
elementary simplifications using \thetag{5}, \thetag{7},
and the braid relations, we get
the following braid monodromy generators:
$$\align
\vec\lambda^+&=\{ A_{4,5,6},\ A_{3,6},\ A_{1,2,6},\
A_{1,3,4},\  A_{2,5}^{A_{3,5}A_{4,5}A_{5,7}},\ A_{4,7},\ A_{1,5,7},\
A_{2,3,7}^{A_{4,7}A_{5,7}A_{3,4}} \},\\
\vec\lambda^-&=\{ A_{4,5,6},\ A_{3,6},\ A_{1,2,6},\
A_{1,3,4},\  A_{2,5},\ A_{4,7}^{A_{5,7}},\ A_{1,5,7},\
A_{2,3,7}^{A_{5,7}}\}.\\
\endalign$$

The braid monodromy presentations of the groups $G^{\pm}=G(\A^\pm)$
may then be found using the Artin representation.  After some
simplifications, we obtain:
$$\split G^+&=\langle u_1,\dots,u_7 \mid [u_4,u_5,u_6], [u_3,u_6],
[u_1,u_2,u_6],[u_1,u_3,u_4], [u_2^{u_3},u_5],\\
&\qquad\qquad\qquad\qquad [u_4,u_7],[u_1,u_5,u_7],[u_2^{u_6},u_3,uv_7]\rangle,\\
G^-&=\langle v_1,\dots,v_7 \mid [v_4,v_5,v_6], [v_3,v_6],
[v_1,v_2,v_6],[v_1,v_3,v_4], [v_2,v_5],\\
&\qquad\qquad\qquad\qquad
[v_4^{v_5},v_7],[v_1,v_5,v_7],[v_2,v_3^{v_5},v_7]\rangle.
\endsplit$$

As mentioned above, $G^+ \cong G^-$.  An explicit isomorphism is given by
$$\matrix
u_1 \mapsto v_4 v_1^{-1}v_4^{-1},\hfill&
u_2 \mapsto (v_6 v_3)^{-1} v_2^{-1} v_6 v_3,\hfill&
u_3 \mapsto v_3^{-1},\hfill&
u_4 \mapsto v_4^{-1}, \hfill\\
u_5 \mapsto v_4 v_5^{-1} v_4^{-1},\hfill&
u_6 \mapsto v_6^{-1}, \hfill&
u_7 \mapsto v_5 v_7^{-1} v_5^{-1}.\hfill
\endmatrix$$
Presentations for $G^{\pm}$ were first obtained by
Rybnikov \cite{Ry}, using Arvola's algorithm.  By Theorem~6.4,
the above presentations are Tietze-I equivalent
to those of Rybnikov.  This can also be seen directly:  For $G^+$, an
isomorphism is given by
$u_1 \mapsto w_7^{-1}$, $u_2 \mapsto w_7 w_4^{-1} w_7^{-1}$,
$u_3 \mapsto w_6^{-1}$, $u_4 \mapsto w_3^{-1}$,
$u_5 \mapsto w_5^{-1}$, $u_6 \mapsto w_1^{-1}$, $u_7 \mapsto w_2^{-1}$,
and similarly for $G^-$.

Since $\A^+$ and $\A^-$ are conjugate, their braid monodromies are
equivalent by Theorem~3.9.  But the two monodromies are not
braid-equivalent.  For, if they were, there would be an isomorphism
$\bar\phi: G^+ \to G^-$ 
determined by a braid automorphism $\phi:F_7\to F_7$ (see Remark~4.3).
In particular,
the induced map on homology, $\bar\phi_*: H_1(G^+) \to H_1(G^-)$,
would be a permutation matrix in $\GL(7,\Z)$, and this is ruled out
by a result of Rybnikov \cite{Ry}, Theorem~3.1.

\enddocument